\documentclass[12pt]{amsart}

\usepackage{amssymb}
\usepackage{amsmath,amsthm,amsfonts, epsfig}
\usepackage{epsfig}
\usepackage{graphicx}
\numberwithin{equation}{section}
\input xy
\xyoption{all}

\def\Title    {Rational points on compactifications of groups}
\def\Author   {Joseph Shalika, Ramin Takloo-Bighash and Yuri Tschinkel}
\def\Subject  {Number Thalika, Ramin Takloo-Bighash and Yuri Tschinkel}
\def\Subject  {Number Theory, Automorphic Forms}
\def\Keywords {Rational points, heights}

\newif\ifpdf
   \ifx\pdfoutput\undefined
   \pdffalse
   \else
           \pdfoutput=1
            \pdftrue
   \fi
      \ifpdf
\usepackage[pdftex,
hyperindex=true,
pdftitle={\Title},pdfauthor={\Author},pdfsubject={\Subject},
pdfkeywords={\Keywords},pdfpagemode={UseOutlines},
bookmarks=true,bookmarksopen=true,
bookmarksopenlevel=0,bookmarksnumbered=true,
pdfstartview={FitV},
colorlinks=true
]{hyperref}
 \else
 \usepackage{hyperref}
 \fi

\usepackage{graphicx}
\advance\headheight 2pt
\calclayout
\allowdisplaybreaks[3]

\theoremstyle{plain}
\newtheorem{thm}{Theorem}[section]
\newtheorem{theo}{Theorem}

\newtheorem{lem}[thm]{Lemma}

\newtheorem{coro}[thm]{Corollary}

\newtheorem{prop}[thm]{Proposition}

\theoremstyle{definition} \theoremstyle{definition}
\newtheorem{defn}[thm]{Definition}

\newtheorem{exam}[thm]{Example}

\newtheorem{rem}[thm]{Remark}           

\newtheorem{notationnum}[thm]{Notation}
\theoremstyle{remark}

\def\vol{{\rm vol}}

\newcommand{\De}{\Delta}
\newcommand{\de}{\delta}
\newcommand{\La}{\Lambda}
\newcommand{\la}{\lambda}
\newcommand{\al}{\alpha}
\newcommand{\C}{\mathbb C}
\newcommand{\Q}{\mathbb Q}
\newcommand{\Z}{\mathbb Z}
\newcommand{\R}{\mathbb R}

\def\cB{{\mathcal B}}
\def\cC{{\mathcal C}}

\def\cH{{\mathcal H}}
\def\cO{{\mathcal O}}
\def\cK{{\mathcal K}}

\def\cL{{\mathcal L}}

\def\cH{{\mathcal H}}

\def\cT{{\mathcal T}}

\def\diag{{\rm diag}}

\def\ca{\check{\alpha}}
\def\Br{{\rm Br}}

\def\sT{{\mathsf T}}
\def\sL{{\mathsf L}}

\def\sn{{\mathsf n}}

\def\sK{{\mathsf K}}
\def\sZ{{\mathsf Z}}
\def\sH{{\mathsf H}}

\def\ga{{\mathbf G}_a}

\def\dv{{\rm div}}

\def\Gr{{\rm Gr}}
\def\Ga{{\Gamma}}
\def\Out{{\rm Out}}
\def\Inn{{\rm Inn}}

\def\G{{\mathsf G}}
\def\H{{\mathsf H}}
\def\M{{\mathsf M}}
\def\L{{\mathsf L}}
\def\U{{\mathsf U}}
\def\P{{\mathsf P}}
\def\T{{\mathsf T}}
\def\S{{\mathsf S}}
\def\B{{\mathsf B}}
\def\D{{\mathsf D}}
\def\V{{\mathsf V}}

\newcommand{\Hom}{{\rm Hom}}
\newcommand{\Ker}{{\rm Ker}}

\newcommand{\Gal}{{\rm Gal}}
\newcommand{\K}{{\mathsf K}}
\def\k{{\mathbf k}}

\newcommand{\PGL}{{\rm PGL}}
\def\Val{{\rm Val}}

\def\Spec{{\rm Spec}}
\def\Lie{{\rm Lie}}
\def\Ad{{\rm Ad}}
\def\pr{{\it pr}}

\def\ra{\rightarrow}

\def\A{{\mathbb A}}
\def\C{{\mathbb C}}

\def\Q{{\mathbb Q}}
\def\rL{{\mathrm L}}

\def\Z{{\mathbb Z}}
\def\C{{\mathbb C}}
\def\N{{\mathbb N}}

\def\ma{{\mathfrak a}}

\def\mg{{\mathfrak g}}
\def\mo{{\mathfrak o}}

\def\mn{{\mathfrak n}}

\def\mO{{\mathfrak O}}

\def\mX{{\mathfrak X}}

\def\fp{{\mathfrak p}}

\def\Pic{{\rm Pic}}
\def\Div{{\rm Div}}

\def\GL{{\rm GL}}

\def\PGL{{\rm PGL}}

\def\zZ{{\mathcal Z}}
\def\eps{{\epsilon}}
\def\tr{{\rm tr}}
\def\X{{\mathfrak X}}
\def\gm{{\mathbb G}_m}

\def\rL{{\rm L}}

\def\rZ{{\mathsf Z}}

\def\bfs{{\mathbf s}}

\def\End{{\rm End}}

\def\lra{\longrightarrow}
\def\ba{\backslash}

\def\sA{{\mathsf A}}

\def\sT{{\mathsf T}}
\def\sL{{\mathsf L}}

\def\sn{{\mathsf n}}
\def\sr{{\mathsf r}}

\def\sK{{\mathsf K}}
\def\sN{{\mathsf N}}

\def\sW{{\mathsf W}}
\def\sw{{\mathsf w}}

\def\sS{{\mathsf S}}
\def\sZ{{\mathsf Z}}

\def\fs{{\mathsf s}}

\def\Nm{{\rm N}}

\makeatother

\makeatletter

\author{Joseph Shalika}
\address{Department of Mathematics\\
Johns Hopkins University\\
3400 N. Charles Str.\\
Baltimore, MD 21218-2686}
\email{shalika@math.jhu.edu}

\author{Ramin Takloo-Bighash}
\address{Department of Mathematics \\
         Princeton University\\
         Fine Hall, Washington Road\\
         Princeton, NJ 08544-1000,  U.S.A.}
\email{rtakloo@math.princeton.edu}
\thanks{The second author was partially supported by the Clay Mathematics
Institute and the NSA}

\author{Yuri Tschinkel}\thanks{The third author was
partially supported by NSF grant 0100277}
\address{Courant Institute, NYU \\
         251 Mercer Str., New York, NY 10012, USA\\
         Mathematisches Institut \\
         Bunsenstr. 3-5 \\
         37073 G\"ottingen, Germany}
\email{tschinkel@cims.nyu.edu}

\title[Rational points]{Rational points on
compactifications of semi-simple groups}

\begin{document}

\date{\today}



\begin{abstract}
We prove Manin's conjecture concerning the distribution of rational points of
bounded height, and its refinement by Peyre, for wonderful compactifications of
semi-simple algebraic groups over number fields. The proof proceeds via the study
of the associated height zeta function and its spectral expansion.
\end{abstract}

\maketitle

\tableofcontents

\setcounter{section}{0}
\section*{Introduction}
\label{sect:introduction}

Let $F$ be a number field and $\G$ a connected linear algebraic
group over $F$. We are interested in distribution properties of
rational points on smooth projective equivariant compactifications
$X$ of $\G$ with respect to heights. More precisely, we assume that
$X$ contains $\G$ as a Zariski open subset and that the left- and
right action of $\G$ on itself extends to $X$. Equivalently, $X$ is
an equivariant compactification of $\G\times \G/\G$. A split over
$F$ semi-simple group $\G$ of adjoint type has a canonical
(wonderful) compactification, constructed over an algebraically
closed field in \cite{concini-p83} and over arbitrary fields (and
over $\Z$) in \cite{strickland} and \cite{deconcini-springer}. In
this paper we consider wonderful compactifications of forms of
such groups. One of our main results is a proof of Manin's
conjecture for this class of varieties:

\begin{theo}
\label{thm:main} Let $X$ be the wonderful compactification of a
semi-simple group $\G$ over $F$ of adjoint type and
$\cL=(L,\|\cdot \|_v)$ an ad\`{e}lically metrized line bundle such
that its class $[L]$ in the Picard group $\Pic(X)$ of $X$ is
contained in the interior of the cone of effective divisors
$\La_{\rm eff}(X)\subset \Pic(X)_{\R}$. Then
\begin{equation}
\label{eqn:zeta} \zZ(\cL,s):=\sum_{x\in \G(F)} H_{\cL}(x)^{-s}=
\frac{c(\cL)}{(s-a(L))^{b(L)}} +\frac{h(s)}{(s-a(L))^{b(L)-1}}
\end{equation}
and, consequently,
$$
\mathcal N(\cL,B):=  \# \{ x\in \G(F)\,|\, H_{{\cL}}(x)\le B\}
\sim \frac{c(\cL)}{a(L)(b(L)-1)!} B^{a(L)} \log(B)^{b(L)-1}
$$
as $B\ra \infty$.
Here
\begin{itemize}
\item $a(L)=\inf\{ a\,\, | \,\, a[L]+[K_X]\in \La_{\rm eff}(X)\}$
(where $K_X$ is the canonical line bundle of $X$);
\item
$b(L)$ is the (maximal) codimension of the face of
$\La_{\rm eff}(X)$ containing $a(L)[L]+[K_X]$;
\item
$c(\cL)\in \R_{>0}$ and
\item
$h(s)$ is a holomorphic function (for $\Re(s)>a(L)-\eps$, some $\eps>0$).
\end{itemize}
Moreover, $c(-\cK_X)$ is the constant defined in \cite{peyre95}.
\end{theo}

\

A concrete special case of this theorem can be described as
follows. Let $\rho: \G \to \GL_N$ be a faithful absolutely
irreducible representation defined over $F$. Since $\G$ is of
adjoint type we may think of $\rho$ as a projective representation
into $\PGL_N$. Think of $\PGL_N$ as an open set in the projective
space $\mathbb{P}^{N^2 - 1}$; fix a height function $H$ on
$\mathbb{P}^{N^2 - 1}(F)$. Denote the pull back of $H$ to $\G(F)$
via $\rho$ by $H_\rho$. The theorem in particular gives an
asymptotic formula for the number of elements of the set
\begin{equation*}
\{ \gamma \in \G(F) \, \, | \, \, H_\rho(\gamma) \leq B \}
\end{equation*}
as $B$ goes to infinity.

The proof is based on the study of analytic properties of the
height zeta function \eqref{eqn:zeta}. The zeta function is first
realized as a special value of an automorphic form on the group
$\G(\A)$. Then we write down the automorphic Fourier expansion in
terms of a basis of automorphic forms. In general, we also have a
contribution from the continuous spectrum. The Fourier expansion
gives an identity of continuous functions which we use to prove a
meromorphic extension of the height zeta function.

This paper is part of a program, initiated in \cite{FMT} and
developed in \cite{batyrev-m90}, \cite{peyre95},
\cite{batyrev-t98}, to relate asymptotics of points of bounded
height to geometric and arithmetic invariants of $X$. For further
results and motivating examples we refer to the book \cite{aster}
and the papers \cite{BT}, \cite{ST}, \cite{chambert-t02} and
\cite{peyre-circle}. The recent preprint \cite{GMO} gives a
simplification of of the proof of Theorem~\ref{thm:main}, based on
ad\`{e}lic mixing.

The paper is organized as follows. In Part One we collect various
results on algebraic groups and automorphic forms.
Section \ref{sect:1-dim} contains a theorem
about transfer of one-dimensional automorphic
representations among inner forms of a group. Section
\ref{spectral-theory} is concerned with Eisenstein series and
automorphic Fourier expansions. The key results
are Lemma~\ref{lem:spectral-expansion-general} and
Proposition~\ref{important-proposition}. In Section~\ref{spherical-functions}
we recall a theorem of Oh (Theorem \ref{Oh})
which we use for estimates in the spectral expansion. In Part
Two we apply this theory to Manin's conjecture.
In Section~\ref{sect:geometry} we review the construction of
wonderful compactifications and their geometric properties.
In Section~\ref{sect:heights} we
study height functions and their
integrals in the local situation. The problem of regularizing global
height integrals is considered in Section
\ref{sect:regularintegral}. In Section~\ref{sect:heightzeta} we
establish analytic properties of height zeta functions in several
complex variables. Tauberian theorems imply asymptotics for the
number of rational points of bounded height.

\

\noindent
{\bf Acknowledgments.}
The authors wish to thank Erez Lapid, Hee Oh, and Peter Sarnak for
various useful discussions.

\part{Algebraic groups and automorphic forms}

\section{Algebraic groups}

\subsection{Basic notation}\label{sect:cartan}

For a number field $F$, let $\Val(F)$ be the set of all places
and $S_{\infty}=S_{\infty}(F)$ the set of archimedean places.
For $v\in \Val(F)$, let $F_v$ be the completion of $F$ with
respect to $v$, $\cO_F$ (resp. $\cO_v)$ the ring of
integers in $F$ (resp. $F_v$) and $q=q_v$ the order of the
residue field $\k=\k_v$ of $\cO_v$.
For any finite set of places
$S$ (containing $S_{\infty}$) we denote by $\cO_S$ the ring of
$S$-integers of $F$. We denote by $\A$ the ring of adeles, by
$\A_S=\prod_{v\notin S}' F_v$ and by $\A_f=\A_{S_{\infty}}$.
Let $\Gamma_F$
and $\Gamma_v$ the Galois groups of $\bar F / F$ and $\bar F_v /
F_v$, respectively. For $L/F$ a Galois extension, with $w
\in \Val(L)$ and $v \in \Val(F)$ satisfying $w | v$,
set $\Gamma_{w/v}=\Gal(L_w/F_v)$.

For $X=X_F$ an algebraic variety over $F$ and $L/F$ an extension
write $X_L=X\times_F L$ for the base-change to $L$
and $X(F)$ for the set of $F$-rational
points of $X$. For $v \in \Val(F)$, we abbreviate $X_v =
X_{F_v}$. We typically identify $X_v(F_v)$ and $X(F_v)$.

We denote by $\Pic(X)$ the Picard group of
$X$ and by $\La_{\rm eff}(X)\subset \Pic(X)_{\R}$ the (closed)
cone of effective divisors on $X$. We often identify
line bundles, divisors and their classes in $\Pic(X)$. If $X$
has an action by an algebraic group $\H$ we write $\Pic^{\H}(X)$
for the group of isomorphism classes of $\H$-linearized line
bundles on $X$.

\subsection{Setup}
\label{setup-intro}
Throughout this paper,
$\G$ is a fixed connected semi-simple group of adjoint type defined
over a number field $F$ and $\G'$ a quasi-split group over $F$ of which $\G$
is an inner form. The groups $\G$ and $\G'$ have a unique split $F$-form
$\G^{sp}$. We may assume that $\G^{sp}$ is obtained by
base-change from a split group, again denoted $\G^{sp}$, over $\Q$.
We fix a Galois extension $E/F$ such that
$\G(E) = \G'(E) = \G^{sp}(E)$. Its Galois
group is denoted by $\Gamma=\Gal(E/F)$. We denote
the simply-connected cover of $\G$ by $\G^{sc}$. Let $\sT$ be a
maximal torus of $\G$ viewed as an algebraic group over $F$. For
each $v\in \Val(F)$, we denote by $\S_v$ the maximal $F_v$-split subtorus
of $\sT$ considered as an algebraic group over $F_v$. The
corresponding objects for $\G'$ and $\G^{sp}$ are written $\sT'$,
$\S_v'$, and $\sT^{sp}$, respectively; for $\G^{sp}$, we make the
choice in such a way that $\sT^{sp}$ is split. We choose
maximal tori $\sT, \sT', \sT^{sp}$ in such a way that they
have the same set of rational points over $E$. In particular,
$\sT'$ is split over $E$. We let $\S'$ be the maximal split torus in
$\G'$ contained in $\sT'$. We choose a finite set $S_E\subset \Val(E)$
containing all the archimedean places so that if $w \notin S_E$
is non-archimedean then
\begin{itemize}
\item
the extension $E_w / F_v$ is unramified for all $v\in \Val(F)$ with $w | v$;
\item
if $v\in \Val(F)$ is not divisible by any $w \in S_E$,
then $\G$ is quasi-split over $F_v$ and $\G(F_v) =
\G'(F_v)$.
\end{itemize}
We will denote by $S_F$ the collection of all places of $F$ which
are divisible by some $w \in S_E$. Later, we will need to enlarge
$S_E$, and $S_F$, to satisfy extra assumptions; c.f.
\ref{sect:goodprimes}. We have a {\em standard} and a {\em
twisted} action of $\Ga=\Gal(E/F)$ on $\G'(E)$, denoted by
$$
g\mapsto \sigma(g), \,\,\, \text{ resp. } \,\,\, g\mapsto
\tilde{\sigma}(g) .
$$
Here
$$
\tilde{\sigma}=c(\sigma)\cdot \sigma,
$$
where $c\in Z^1(\Gamma, \Inn(\G')(E))$ is a cocycle with values in
the group of inner automorphisms of $\G'$.

Pick a Borel subgroup $B^{sp}$ in $\G^{sp}$ which contains
$T^{sp}$. As $\G'$ is quasi-split it contains a Borel subgroup
$\B'$, containing $\sT'$, which is defined over $F$. If $v \notin
S_F$, then by transfer of structure we obtain a Borel subgroup
$\B(F_v)$ in $\G(F_v)$. If $v \in S_F$, we pick a minimal
parabolic $F_v$-subgroup $\B(F_v)$ of $\G(F_v)$ containing
$\S(F_v)$. For each place $v$, denote by $\tilde\Phi(\G(F_v),
\S(F_v))$ the set of roots of $\S(F_v)$ in $\G(F_v)$ and by
$\Phi(\G(F_v), \S(F_v))$ the set of non-multipliable roots in
$\tilde\Phi$ with the ordering given by $\B(F_v)$. Let
$\mathfrak{X}^*(\S(F_v))$ denote the set of characters of
$\S(F_v)$ defined over $F_v$. Denote by $\mathfrak{X}^+$ (resp.
$\Phi^+$) the set of positive characters (resp. roots) in
$\mathfrak{X}^*(\S(F_v))$ (resp. $\Phi(\G(F_v), \S(F_v))$) with
respect to that ordering. We also let $\Delta(\G^{sp}, \sT^{sp}) =
\{\alpha_1, \dots, \alpha_r\}$ be the set of simple roots of for
the pair $(\G^{sp}, \sT^{sp})$ with respect to the ordering
introduced by $\B^{sp}$, and let $\rho$ be half the sum of the
positive roots for the same pair. For any subset $I \subset [1,
\dots, r]$ we denote by $\P_I^{sp}$ the corresponding parabolic
subgroup of $\G^{sp}$ and by $\L_I^{sp}$ its Levi subgroup. Let
$\sW$ be the Weyl group of the pair $(\G^{sp}, \sT^{sp})$. Denote
by $\{\fs_1,...,\fs_\sr\}$ the set of simple reflections, by
$\ell$ the length function on $\sW$ and by $\sw_0\in \sW$ the
longest element with respect to $\ell$. We let $\{\omega_1, \dots,
\omega_r\}$ be the collection of the fundamental weights of the
simply-connected cover of $\G^{sp}$. We will also occasionally
speak of the root system of the pair $(\G', \sT')$, the ordering
given by the choice of $\B'$, and the collection of simple roots
$\Delta(\G', \sT')$; the collection of $\Gamma_F$-orbits in
$\Delta(\G', \sT')$ is the same as $\Delta(\G', \S')$.

If $v$ is archimedean, we set
\begin{equation}
F_v^0 = \{ x \in \R \, | \, x \geq 0 \}
\end{equation}
and
\begin{equation}
\hat{F}_v = \{ x \in \R \, | \, x \geq 1 \}.
\end{equation}
If $v$ is non-archimedean, we fix a uniformizer $\varpi$ of $F_v$,
and set
\begin{equation}
F_v^0 = \{ \varpi^n \, | \, n \in \Z \}
\end{equation}
and
\begin{equation}
\hat{F}_v = \{ \varpi^{-n} \, | \, n \in \N \}.
\end{equation}
We set
\begin{equation}
\S_v(F_v)^0 = \{ a \in \S_v(F_v) \, | \, \alpha(a) \in F_v^0
\text{ for each } \alpha \in \mathfrak X^*(\S(F_v))\}
\end{equation}
and
\begin{equation}
\S_v(F_v)^+ = \{ a \in \S(F_v) \, | \, \alpha(a) \in \hat{F}_v
\text{ for each } \alpha \in \Phi^+ \}.
\end{equation}
Similarly we may define $\S_v'(F_v)^+$. We define constants
$\kappa_\alpha$ for $\alpha \in \Delta(\G^{sp}, \T^{sp})$ by:
\begin{equation}
\sum_{\al>0,\al\in \Phi(\G,\T)}\al =\sum_{\al\in \Delta(\G,\T)}
\kappa_{\al} \al.
\end{equation}

\subsection{Integrality}
\label{sect:local}

By assumption, $\G$ and $\G'$ are of adjoint type and we can
identify $\Inn(\G')=\G'$. Then, for $g\in \G'(E)$, we have
$$
c(\sigma)\cdot g=a_{\sigma} ga_{\sigma}^{-1}
$$
for a uniquely defined element $a_{\sigma}\in\G'(E)$.
%
For $w \in \Val(E)$ and $v \in \Val(F)$ with $w | v$,  we have a
restriction homomorphism
$$
\begin{array}{cccc}
\iota_{w} \,:\, &  Z^1(\Ga,\G'(E)) & \ra  & Z^1(\Ga_{w/v},\G'(E_w)), \\
                &     c                & \mapsto & c_w.
\end{array}
$$
The group $\G_v=\G_{F_v}$ is obtained from the corresponding
quasi-split group $\G'_v$ by twisting with $c_w$. With our choice
of the set $S_E$ from \ref{sect:cartan}, the local cocycle $c_w$
splits for $w\notin S_E$. Thus for $w\notin S_E$, $\sigma\in
\Ga_w$, we have
$$
c_{w}(\sigma) = a_w \cdot \sigma(a_{w}^{-1}) \ ,
$$
for some element $a_w\in \G'(E_w)$. Identifying $\G'(E)=\G(E)$ we
may regard
\begin{itemize}
\item $\G'(F)$ as the $\Ga$-fixed points of $\G'(E)$ for the
standard action and \item $\G(F)$ as the $\Ga$-fixed points of
$\G'(E)$ for the twisted action.
\end{itemize}
Locally, we have a similar situation:

\

\centerline{ \xymatrix{
\G(E_w) \ar@{=}[r] & \G'(E_w) \\
\G(F_v) \ar@{->}[u] & \G'(F_v) \ar@{->}[u]
 }
}

\

\noindent and further
$$
a_w^{-1}\G(F_v)a_w =\G'(F_v)\ .
$$

Recall that the global group $\G'$ splits over $E$ and put
$$
\mg'=\Lie(\G') \,\, \text{ and }\,\, \mg'_E=\mg'(E) \ .
$$
Let $\T'\subset \G'$ be a maximal split torus defined over $E$.
Let
$$
\mg^{sp}=\Lie(\G^{sp}), \,\,\, \mg_{\Q}^{sp}=\mg^{sp}(\Q) \ .
$$
Let $\La_{\Q}$ be a Chevalley lattice in $\mg_{\Q}^{sp}$ adapted to
$\T^{sp}$. Set
$$
\La_E':=\La_{\Q} \otimes_{\Z} \cO_E
$$
and also set
$$
\La_w':=\La_{E}'\otimes{\mo_E} \cO_w.
$$
Then $\La_w'$ is a Chevalley lattice in $\mg_w'=\mg'(E_w)$.
Moreover, $\La_w'$ is adapted to the torus $\T'_w$ obtained from
$\T'$ by base extension to $E_w$.
Let $\K_w$ be the stabilizer of $\La_w$ in
$\G(E_w)=\G'(E_w)$. We have (see \cite{im}, Proposition 2.33):
$$
\G_w=\K_w\cdot \B(E_w)=\B(E_w)\cdot \K_w \ .
$$

\begin{prop}
\label{prop:11} Let $v$ be a finite place of $F$ and $E_w$ an
unramified extension of $F_v$. An element $a_w\in
\G(E_w)=\G'(E_w)$ satisfies
$$
a_w\sigma(a_w)^{-1} \in \K_w
$$
if and only if it is of the form
$$
a_w=k_w\cdot \gamma_w,
$$
with $k_w\in \K_w $ and $\gamma_w\in \G'(F_v)$.
\end{prop}

We first recall the following

\begin{thm}[see \cite{pr}, p. 292]
\label{theo:pr} Let $\G_v$ be a connected semi-simple group
over $F_v$. Suppose that $\G_v$ is defined over
$\mo_v$ and that as a group scheme over $\mo_v$, it
has a connected smooth reduction modulo $\mathfrak p_v$. Let $E_w$
be a finite unramified extension of $F_v$. Then
$$
H^1(\Ga_w,\G(\cO_w))=1.
$$
\end{thm}

\begin{proof}[Proof of Proposition~\ref{prop:11}]
Consider the exact sequence of (pointed) sets
$$
1\ra \G_v'(\cO_w)\ra \G_v'(E_w)\ra \G_v'(\cO_w)\ba \G_v'(E_w)\ra 1
$$
and the corresponding exact sequence in cohomology
$$
1\ra \G_v'(\cO_v)\ra \G_v'(F_v)\ra (\G_v'(\cO_w)\ba
\G_v'(E_w))^{\Ga_w}\ra H^1(\Ga,\G_v'(\cO_w))=1,
$$
by Theorem~\ref{theo:pr}. Consequently,
$$
\G_v'(\cO_v)\ba \G_v'(F_v) = \left(\G_v'(\cO_w)\ba
\G_v'(E_w)\right)^{\Ga_w} \ ,
$$
precisely the statement of Proposition~\ref{prop:11}.
\end{proof}

\subsection{Good primes}
\label{sect:goodprimes} To apply Proposition~\ref{prop:11} we will
enlarge $S_E$, and $S_F$, so that for all $v\notin S_F$ one has:
\begin{itemize}
\item $\G'$ has a group-scheme structure via the Chevalley lattice
$\La_F$; \item $\G'$ has a connected smooth reduction modulo
$\fp_v$; \item for all $\sigma \in \Ga$
$$
c(\sigma)\in \G'(\cO_w)=\K'_w \ ;
$$
\item the local cocycle
$$
c_w\in Z^1(\Ga_w,\G'(E_w))
$$
splits.
\end{itemize}

Under these assumptions, for $\sigma\in \Ga_w$ we may write
$$
c_w(\sigma)=a_w\cdot \sigma(a_w^{-1})
$$
for some $a_w\in \G'(E_w)$. Then, for $w\notin S$,
$$
a_w\cdot\sigma(a_w^{-1})\in  \K'_w \ .
$$
By Proposition~\ref{prop:11}, we have $ a_w=k_w\cdot \gamma_w, $
for some $k_w\in \K'_w, \gamma_w\in \G'(F_v)$. Since
$\sigma(\gamma_w)=\gamma_w$ for $\sigma\in \Ga_w$, we may assume
that
$$
c_w(\sigma)=a_w\cdot \sigma(a_w^{-1})
$$
with $a_w\in \K'_w$. Our main conclusion is

\begin{coro}
\label{coro:comp} There is a unique element $a_w\in \K'_w$ such
that
$$
a_w^{-1}\G(F_v)a_w=\G'(F_v).
$$
\end{coro}
\begin{proof}
We only need to prove uniqueness. Suppose that we have another
element $b\in \G(E_w)$ such that
$$
b^{-1}\G(F_v)b=\G'(F_v).
$$
Let $z=ba^{-1}$. Then
\begin{equation}
\label{eqn:1.15.1} z\G(F_v)z^{-1} =\G(F_v).
\end{equation}
Next, let $\sigma\in \Gal(E_w/F_v)$. Then from \eqref{eqn:1.15.1}
we have
$$
\tilde\sigma(z g z^{-1}) = zgz^{-1}
$$
for all $g\in \G(F_v)$. Hence, $z^{-1}\sigma(z) $ belongs to the
centralizer of $\G(F_v)$ in $\G(\bar{F}_v)$. But then by a theorem
of Rosenlicht, $\G(F_v)$ is Zariski dense in $\G(\bar{F}_v)$ (see
\cite{borel}, p. 11). Thus $z\in \rZ(\bar{F}_v)$ (the center of
$\G(\bar{F}_v)$). Since $\G$ is an adjoint group, we have then
$z=1$ and $a=b$.
\end{proof}
We set $\K'_v = \K'_w \cap \G'(F_v)$. For $v \notin S_F$, this is
a maximal compact subgroup in $\G'(F_v)$. Put $\K_v = a_w \K'_v
a_w^{-1}$. Then
\begin{coro}[Cartan Decomposition]
For $v \notin S_F$, we have
\begin{equation}
\G(F_v) = \K_v \S(F_v)^+ \K_v.
\end{equation}
\end{coro}
\begin{rem}
If $v \in S_F$, then there is a maximal compact subgroup $\K_v$ of
$\G(F_v)$ and a finite set $\Omega_v \subset \G(F_v)$ such that
\begin{equation}
\G(F_v) = \K_v \S(F_v)^+ \Omega_v \K_v.
\end{equation}
The set $\Omega_v$ is trivial if $v$ is archimedean, or $\G(F_v)$
is split, or $\G(F_v)$ quasi-split and split over an unramified
extension.
\end{rem}

\subsection{Galois cohomology of simply-connected groups}
\label{sect:groups-adele}

Here we recall basics facts about simply-connected semi-simple
groups over the adeles.

\begin{thm}[Kneser]
\label{thm:kneser} Let $\H$ be a simply connected semi-simple
group over $F$. Then, for $v\nmid \infty$,
$$
H^1(\Gamma_v,\H(\bar{F}_v))=1
$$
and the natural map
$$
H^1(\Gamma, \H^{sc}(F))\ra \prod_{v|\infty} H^1(\Gamma_v,
H(\bar{F}_v)
$$
is a bijection.
\end{thm}


It is useful to keep in mind the following:

\begin{lem}[Theorem 6.17 of \cite{pr}]
\label{lemm:usefl} Let $\U$ be a connected compact algebraic
group over $\R$. Then
$$
H^1(\Gamma_{\R}, \U(\C))=\sA'/\mathsf W_\sA
$$
where
$$
\sA' :=\{ t\in \sA\, |\, t^2 =1\},
$$
$\sA\subset \U$ is the maximal $\R$-torus and $\mathsf
W_\sA=\mathsf N_\U (\sA)/\sA$ is the corresponding Weyl group.
Here, $N_\U(\sA)$ is the normalizer of $\sA$ in $\U$.
\end{lem}

\subsection{Norm maps}
\label{sect:norms}

Let $\T$ be an algebraic torus over $F$, splitting over a normal
extension $E/F$ with Galois group $\Gamma$, $v$ a place of $F$,
$w\in \Val(E)$ a place over $v$, and $\Gamma_{w/v}$ the
corresponding local Galois group.
We have the following natural {\em norm } homomorphisms:
$$
\begin{array}{cclcccl}
\T(E_w)  & \stackrel{\Nm_{w/v}}{\longrightarrow} & \T(F_v), &  &
\T(E_w)  &\stackrel{\Nm_{v}}{\longrightarrow}    & \T(F_v)     \\
 t_w     & \mapsto                               & \prod_{\sigma\in \Gamma_{w/v} } \sigma(t_w),  & &
 t_w     & \mapsto                               & \prod_{w\mid v} \Nm_{w/v}(t_w), \\
\T(E)    & \stackrel{\Nm_{E/F}}{\longrightarrow} & \T(F),  &   &
\T(\A_E) & \stackrel{\Nm_{\A_E/\A_F}}{\longrightarrow} & \T(\A_F) \\
  t      & \mapsto                               & \prod_{\sigma\in \Gamma} \sigma(t), &  &
(t_w)_w  & \mapsto & (\Nm_v(t_w))_v.
\end{array}
$$

\begin{lem}
\label{lemm:normcomp} For $t\in \T(E)$ one has $ \Nm_{E/F} (t)
=\Nm_{\A_E/\A_F} (t). $
\end{lem}

\begin{proof}
Let $w$ be a place of $E$ over $v$ and $\Ga_w\subset \Gamma$ the
fixer of $w$. For any $\sigma \in \Gamma$ we have
$\sigma\Ga_w\sigma^{-1} =\Ga_{\sigma(w)}$ and $ |x|_w
=|\sigma(x)|_{\sigma(w)}. $ For a fixed $w$ write
$$
\Ga =\cup_{j=1}^g \sigma_j \Ga_{w_1} =\cup_{j=1}^g \Ga_{w_1}
\sigma_j^{-1}.
$$
We have, for $t\in \T(E)$,
\begin{align*}
\Nm_{E/F} (t) = & \prod_{j=1}^g \prod_{\rho\in \Ga_{w_1}} \rho\sigma_j^{-1}(t)\\
              = & \prod_{j=1}^g\sigma_j^{-1} \prod_{\rho\in
                \Ga_{w_1}} \sigma_j \rho\sigma_j^{-1}(t) \\
              = & \prod_{j=1}^g \sigma_j^{-1} \prod_{\rho\in \Ga_{w_j}} \rho(t).
\end{align*}
We have the commutative diagram

\

\centerline{ \xymatrix{ E \ar[r]_{\iota_w}  \ar[d]^{\sigma} &
E_w\ar[d]^{\sigma^*}  &
 F_v\ar[l]\ar[d]^{id} & \ar[l]^{\iota_v}F \\
  E \ar[r]_{\iota_{\sigma(w)}}        & E_{\sigma(w)} &
 F_v  \ar[l]  &
} }

\

\noindent where $\iota_v,\iota_w$ are the canonical maps and
$\sigma^*$ is the continuous extension of $\sigma$.

For $t_w:=\prod_{\rho\in \Ga_{w}} \rho(t)$ we have
\begin{align*}
\iota_v(\Nm_{E/F} (t)) & = \iota_{w_1} (\Nm_{E/F}(t)) \\
                       & = \prod_{j=1}^g \iota_{w_1}\sigma_j^{-1} (a_j) \\
                       & = \prod_{j=1}^g (\sigma_j^*)^{-1}\cdot \iota_{w_j}(a_j).
\end{align*}
Next, let for each $w$,
$$
\Ga_w^*:=\{ \rho^*\,|\, \rho\in \Ga_w\}.
$$
If $\rho\in \Ga_{w_j}$ then $\rho(a_j)=a_j$. Thus
$$
\iota_{w_j}(a_j)=\iota_{w_j}\rho(a_j) =\rho^*(\iota_{w_j}(a_j)).
$$
Thus $\iota_{w_j}(a)$ is fixed by $\Ga_{w}^*$ and therefore
belongs to $F_v$. But then
$$
(\sigma_j^*)^{-1}(\iota_{w_j}(a_j)) = \iota_{w_j}(a_j).
$$
Thus
$$
\iota_v(\Nm_{E/F}(t))=\prod_{j=1}^g \iota_{w_j}(a_j).
$$
On the other hand,
\begin{align*}
\Nm_{\A_E/\A_F} (t) & =   \prod_{j=1}^g \Nm_{E_{w_j}/F_v}(\iota_{w_j}(t))\\
                    & = \prod_{j=1}^g \prod_{\rho^*\in \Ga_{w_j}^*} \rho^*\iota_{w_j}(t)  \\
                    & = \prod_{j=1}^g  \prod_{\rho\in \Ga_{w_j}} \iota_{w_j} (\rho(t))  \\
                    & = \prod_{j=1}^g \iota_{w_j} (a_j).
\end{align*}
This completes the proof.
\end{proof}

\section{One-dimensional automorphic representations}
\label{sect:1-dim}

In this section we define a transfer associating to
each one-dimensional automorphic  representation $\chi $ of $\G(\A)$
a one-dimensional automorphic representation $\chi'$ of $\G'(\A)$ in
such a way that it respects local isomorphisms. The basic idea is to
define a local transfer, and then use invariance under $\G(F)$ and
weak approximation to extend the definition to the entire
$\G(\A)$. We will also define certain Euler products that will later
be used to regularize height integrals.

\subsection{Automorphic characters for simply-connected groups}
\label{sect:chi-sc}

\begin{prop}
\label{prop:1} Let $F$ be a number field and $\H$ a connected
simply-connected semi-simple group over $F$. Let $\chi$ be a
one-dimensional automorphic representation of $\H(\A)$. Then
$\chi=1$.
\end{prop}

\begin{proof}
Let $E/F$ be a finite Galois extension splitting $\H$.
There are infinitely many  $v\in \Val(F)$ which split completely in $E$.
For such $v$
the local group $\H(F_v)$ is split. Since $\H$ is
simply-connected, $\H(F_v)$, as an abstract group, is its own
derived group (see \cite{borel-alg}, esp. 3.3.5, p. A-16). Thus
$\chi|_{\H(F_v)}=1$. Also, $\chi|_{\H(F)}=1$. By strong
approximation (see Theorem 7.12, p. 427 of \cite{pr}),
$\H(F_v)\H(F)$ is dense in $\H(\A)$. It follows that $\chi=1$.
\end{proof}

\subsection{The local transfer}
Let $\G$ be as in \ref{setup-intro} and
$$
j_v\,:\, \G^{sc}(F_v)\ra \G(F_v)
$$
be the canonical homomorphism from its simply-connected covering.
Let $v$ be a finite place of $F$ and $\hat{\G}_v$ the group of
characters $\chi_v$ of $\G(F_v)$ which are trivial on
$j_v(\G^{sc}(F_v))$.

By Kneser's theorem \ref{thm:kneser}, $H^1(F_v,
\G^{sc}(\bar{F}_v))=1$.
We have an exact sequence
$$
1\ra \rZ(F_v)\ra \G^{sc}(F_v)\stackrel{j_v}{\longrightarrow}
\G(F_v) \stackrel{\delta_v}{\longrightarrow} H^1(F_v,
\rZ(\bar{F}_v))\ra 1
$$
and a similar sequence for $\G'$. Here $\rZ$ is the center of
$\G^{sc}$.
We may identify
\begin{equation}
\label{eqn:1.5.2} \hat{\G}_v = H^1(F_v,\rZ(\bar{F}_v))^{\vee},
\end{equation}
(the character group of $H^1(F_v, \rZ(\bar{F}_v))$. Recall that
$\G$ is obtained from $\G'$ by replacing the standard Galois
action by a twisted action, via some representative of
$H^1(\Gal(\bar{F}/F), \Inn(\G'(\bar{F})))$. Since inner
automorphisms of $\G'(\bar{F})$ fix $\rZ(\bar{F})$, these two
actions coincide on $\rZ(F)$, and similarly for $\rZ(\bar{F}_v)$.
By \eqref{eqn:1.5.2} we get then, for non-archimedean $v$, a
natural isomorphism
\begin{equation}
\label{eqn:tr} \tr_v\,:\, \hat{\G}_v  \ra \hat{\G}_v'.
\end{equation}
\begin{rem}
We observe that equation \eqref{eqn:1.5.2} implies that there is a
number $n$, independent of $v$, such that every occurring
character of $\G(F_v)$ satisfies $\chi^n =1$. This follows from a
case-by-case analysis of $H^1(F_v, Z(\overline{F}_v))$ as carried
out in \cite{pr}, Section 6.5. We will come back to this
point in \ref{ex:pgln}.
\end{rem}

\subsection{Weak approximation}

We now apply weak approximation (Theorem 7.7,
p. 415 of \cite{pr}):

\begin{thm}
The group $\G(F)$ is dense in $\G_{\infty}:=\prod_{v\mid \infty}
\G(F_v) $.
\end{thm}

For $v |\infty$, the map $j_v\,:\, \G^{sc}(F_v)\ra \G(F_v)$ is
submersive. Write $\G_{\infty}^*\subset \G_{\infty}$ for the image
of $\prod_{v |\infty} \G^{sc}(F_v)$ under $\prod_{v\mid \infty}
j_v$, it is an open subgroup. Thus $\pr_{\infty}(\G(F))\cdot
\G_{\infty}^*= \G_{\infty}$ and
$$
\G(\A)=\G(F)\cdot \G_0\cdot \G_{\infty}^*.
$$
Here $\pr_\infty$ is the product of the projection maps $\pr_v:
\G(F) \to \G(F_v)$ for $v | \infty$.

\subsection{Global transfer} Let $\chi=\prod_v \chi_v$ be a one-dimensional
automorphic representation of $\G(\A)$ such that $\chi_v$ is
trivial on $j_v(\G^{sc}(F_v))$, for all $v\in \Val(F)$. In
particular, $\chi_{\infty}: = \prod_{v | \infty} \chi_v$ is
trivial on $\G_{\infty}^*$. Using \eqref{eqn:tr} we define
$$
\chi_0' =\prod_{v\nmid \infty} \tr_v(\chi_v).
$$
We also have
$$
\G'(\A) =\G'(F)\cdot \G_0' \cdot (\G_{\infty}')^*.
$$
We extend $\chi_0'$ to a character $\chi'$ of $\G'(\A)$ by setting
$\chi'=1$ on $\G'(F)\cdot (\G_{\infty}')^*$. Then $\chi'$ is
well-defined  since it is trivial on the intersection
$$
\G_{0}'\cap \G'(F)\cdot (\G'_{\infty})^*.
$$

\

Indeed, write $ \G(F)^*:=\{ \gamma\in \G(F)\,|\,
\pr_\infty(\gamma) \in \G_{\infty}^*\}, $ and similarly for $\G'$.
Note that
$$
\pr_0(\G(F)^*)=\G(\A_f)\cap \G(F)\cdot \G_{\infty}^*,
$$
and similarly for $\G'$. It suffices to check that
\begin{equation}
\label{eqn:1.8.1}
\delta_0(\pr_0(\G(F)^*)=\delta_0'(\pr_0(\G'(F)^*)),
\end{equation}
where $\delta_0:=\prod_{v\nmid \infty} \delta_v$ and similarly,
$\delta_0'=\prod_{v\nmid \infty}\delta_v'$. In fact, the character
$\chi$ is trivial on $\G(F)\cdot \G_{\infty}^*$.

Next, the sequence
\begin{equation}
\label{eqn:1.8.2} \G^{sc}(F_v)\stackrel{j_v}{\longrightarrow}
\G(F_v)\stackrel{\delta_v}{\lra} H^1(F_v, \rZ(\bar{F}_v))
\end{equation}
is exact for all $v$. Set $\delta_{\infty}=\prod_{v|
\infty}\delta_v$. Then in particular,
$\delta_{\infty}(\G_{\infty}^*)=\{ e\} $. We also set $\delta_{\A}
= \delta_0\times \delta_{\infty} = \prod_{v} \delta_v$. Then
\eqref{eqn:1.8.1} is equivalent to
\begin{equation}
\label{eqn:1.9.1} \delta_{\A}(\G(F)^*) =\delta_{\A}'(\G'(F)^*).
\end{equation}

Let $p_v$ be the canonical map
$$
p_v\,:\, H^1(F,\rZ(\bar{F}))\ra H^1(F_v, \rZ(\bar{F}_v))
$$
and
$$
p\,:\, H^1(F,\rZ(\bar{F}))\ra \prod_{v}H^1(F_v,\rZ(\bar{F}_v))
$$
the corresponding global map. We have $p\circ \delta_F
=\delta_{\A}$, where $\delta_F$ is the coboundary map $\G(F) \ra
H^1(F,\rZ(\bar{F}))$.

\begin{prop}
Let $c\in H^1(F,\rZ(\bar{F}))$. Then $c\in \delta_F(\G(F)^*)$ if
and only if
$$
p_v(c)=c_v=0\,\,\, \text{ for  all } \,\,\, v\mid \infty.
$$
\end{prop}

\begin{proof}
We have a diagram

\

\centerline{ \xymatrix{
  \G(F) \ar[r]^<<<<<{\delta_F}  \ar[d] &  H^1(F,\rZ(\bar{F}))\ar[d]_{z}\ar[r]^<<<<<<<<<{k_F}&
H^1(F,\G^{sc}(\bar{F}))\ar[d]^{g}\\
  \G(\A)  \ar[r]^<<<<<{\delta_{\A}}  &  \prod_v H^1(F_v,\rZ(\bar{F}_v)) \ar[r]^<<<<{k_{\A}}&
\prod_{v\nmid \infty} H^1(F_v, \G^{sc}(\bar{F}_v)). } }

\

 \noindent
where $g$ is bijective (see Th. 6.6, p. 286 of \cite{pr}), both
squares are commutative and the top row is exact. Recall also that
$H^1(F_v, \G^{sc}(\bar{F}_v))=\{ e\}$,  if $v$ is finite.

Suppose that $c_v=z_v(c)=0$ for $v$ infinite, $c\in
H^1(F,\rZ(\bar{F})$. Then $k_{\A}((c_v)_v)=e$. Consequently, $k_{\A}
\cdot z (c) = g \cdot k_F(c) = e$. Since $g$ is bijective, we have
$k_F(c)=0$ and consequently
$$
c=\delta_F(\gamma)
$$
for some $\gamma\in \G(F)$. We prove that $\gamma\in \G(F)^*$. In
fact, we have
$$
c_v=\delta_v(\gamma_v)=e, \,\,\, \text{ for }\,\,\, v\mid \infty.
$$
But then, by the exactness of the sequence \eqref{eqn:1.8.2},
$\gamma_v\in {\rm Im}(j_v)$, for any infinite $v$. Therefore,
$\gamma\in \G^*_{\infty}$. The converse assertion also follows
from the exactness of \eqref{eqn:1.8.2}.
\end{proof}

\begin{coro}
As before, let
$$
\delta_F\,:\, \G(F)\ra H^1(F,\rZ(\bar{F}))
$$
and similarly $\delta_F'$ be the respective co-boundary maps. Then
$$
\delta_F(\G(F)^*) =\delta_F'(\G'(F)^*).
$$
\end{coro}

\begin{proof}
By the Proposition, each of the sets $\delta_F(\G(F)^*)$ and
$\delta_F'(\G(F)^*)$ coincides with the set of $c\in
H^1(F,\rZ(\bar{F}))$ satisfying $c_v=e$ for $v\mid \infty$. Our
assertion is now an immediate consequence. \end{proof}

\subsection{Compatibility} The following lemma shows that the above definition of
transfer is the correct one:
\begin{lem}
\label{lemm:chii} Let $v \notin S_F$ and $w \in \Val (E)$ be such
that $w\mid v$. Let $a_w$ be as in Corollary \ref{coro:comp}. Then for
all $x\in \G(F_v)$ one has
$$
\chi_v'(a_w^{-1}xa_w)=\chi_v(x).
$$
\end{lem}

\begin{proof}
In this proof, we write $a$ for $a_w$. There is a character $\xi$
of the discrete group $H^1(F_v, \rZ(\bar{F}_v))$ such that
$$
\chi(x)=\xi(\delta_v(x)), \,\,\, x\in \G(F_v) \,\,\, \text{ and
}\,\,\, \chi'(x)=\xi(\delta_v'(x)), \,\,\, x\in \G'(F_v).
$$
It suffices to prove that $ \delta_v'(a^{-1}xa)=\delta_v(x). $ We
calculate the two co-boundary maps, $\delta_v$ and $\delta_v'$.
Choose $\tilde{x}\in \G^{sc}(\bar{F}_v)$ such that
$j_v(\tilde{x})=x$. Then
$$
\delta_v(x) =\tilde{x}^{-1}\sigma(\tilde x)=
\tilde{x}^{-1}c(\sigma)\circ \sigma'(\tilde x) =
\tilde{x}^{-1}\tilde{a}\sigma'(\tilde a)^{-1} \sigma'(\tilde{x})
\sigma'(\tilde{a}) \tilde{a}^{-1}.
$$
Recall that, for $v\notin S_F$,
$$
c_v(\sigma)=a^{-1}\sigma'(a).
$$
We calculate $\delta_v'(a^{-1}xa)$, for $x\in \G(F_v)$. Choose
$\tilde{a}\in \G^{sc}(\bar{F}_v)$ such that $\Ad(\tilde{a})=a$. We
have then
$$
\delta_v' (a^{-1} xa) = (\tilde{a}^{\-1} \tilde x\tilde a)^{-1}
\sigma'(\tilde{a}^{-1} \tilde x\tilde{a})= \tilde{a}^{\-1} \tilde
{x}^{-1}\tilde a \sigma'(\tilde{a}^{-1})
\sigma'(\tilde{x})\sigma'(\tilde{a}).
$$
Since $\delta_v'(x)\in \rZ(\bar{F}_v)$, we have
$$
a\delta_v'(x) =\delta_v'(x)a,
$$
or
$$
\tilde{a}\delta_v'(x)\tilde{a}^{-1} =\delta_v'(x).
$$
Thus
$$
\delta_v' (a^{-1} xa) =\tilde{x}^{-1}
\tilde{a}\sigma'(\tilde{a}^{-1})
\sigma'(\tilde{x}^{-1})\sigma'(\tilde{a})\tilde{a}^{-1} =
\delta_v(x),
$$
as claimed.
\end{proof}

\subsection{Continuity}

Since $\rZ$ is commutative, $H^1(F_v, \rZ(\bar{F}_v))$ is a
(discrete) group and the co-boundary map
$$
\delta_v\,:\, \G(F_v)\ra H^1(F_v, \rZ(\bar{F}_v))
$$
is a group homomorphism.

\begin{lem}
The map $\delta_v$ is continuous.
\end{lem}

\begin{proof}
The map
$$
j_v\,:\, \G^{sc}(F_v)\ra \G(F_v)
$$
is a submersion of $v$-adic Lie groups
(the differential of $j_v$ is the identity on the Lie algebra
$\Lie(\G^{sc}(F_v)=\Lie(\G(F_v))$. Hence, ${\rm
Im}(j_v)=\Ker(\delta_v)$ is an open subgroup of $\G(F_v)$ and the
claim follows.
\end{proof}

In any case, if $\chi_v\,:\, \G(F_v)\ra \C^*$ is a homomorphism of
groups which is trivial on the image of $j_v$, then $\chi_v$ is
automatically continuous.

\subsection{Automorphy}
We have associated to each one-dimensional automorphic
representation $\chi$ of $\G(\A)$, trivial on $j(\G^{sc}(\A))$, a
``formal'' automorphic character
$$
\chi'\,:\, \G'(\A)\ra \C^\times,
$$
i.e., a homomorphism, trivial on $\G'(F)$.

\begin{thm}
The character $\chi'$ is an automorphic character.
\end{thm}
\begin{proof}
We need to show that $\chi'$ is continuous. Since the map
$$
\delta_v'\,:\, \G'(F_v)\ra H^1(F_v,\rZ(\bar{F}_v))
$$
is continuous, each $\chi_v'$ is also continuous, for $v$ finite.
By construction, $\chi_{\infty}'=\prod_{v\mid \infty} \chi_v'$ is
trivial on the open subgroup
$(\G'_{\infty})^*\subset\G_{\infty}'$. Thus, $\chi_{\infty}'$ is
continuous. The character $\chi_v$ is trivial on $\K_v$, for
almost all $v$. We claim that $\chi'$ is trivial on
$\K_v'$ for almost all $v$. This immediately follows from Lemma
\ref{lemm:chii} and considerations of Section~\ref{sect:goodprimes}.
We have already noted that $\chi'_v$ is continuous for all $v$. It
follows finally that $\chi'$ is continuous.
\end{proof}

\subsection{Hecke characters}
\label{sect:1.6}
We now associate to $\chi'$ a function $L(s,\chi')$,
a product of Hecke $L$-functions. This function will later be
used to regularize Fourier transforms of global height functions.

Suppose first that $\G'_F = \G^{sp}$ is split. Let
$\T^{sp}\subset \B^{sp}\subset \G^{sp}$ be a maximal split torus
(over $F$), and a Borel subgroup containing $\T$. Let $\Delta =
\{\alpha \}$ be the associated set of simple roots and $\{\ca\}$
the dual basis of the co-characters $\X_*(\T^{sp})$. For $\la\in \gm(\A)$ and
$\al\in \Delta$ define
$$
\xi_{\al}(\la) = \chi'(\ca(\la)).
$$
This $\xi_\al$ is a character of $\gm(\A)/\gm(F)$
(Hecke-character). We define
$$
L(s,\chi') =\prod_{\al\in \Delta} L(s,\xi_\al).
$$
In general, when $\G'$ is quasi-split, and not necessarily split
over $F$, we have
$$
\Phi(\G',\T')\subset \X^*_E(\T')=\X^*(\T')\,\, \text{ and } \,\,
\Phi(\G',\sS')\subset \X^*_F(\T').
$$
Then the restriction map
\begin{equation}
\label{eqn:1.31.1} r\,:\, \Phi(\G',\T')\ra \Phi(\G',\sS')
\end{equation}
is surjective, and
$$
r(\Delta(\G',\T'))=\Delta(\G',\sS').
$$
The Galois group $\Ga$ acts transitively on the fibers of
$$
r\,:\, \Delta(\G',\T')\ra \Delta(\G',\sS').
$$
Fix a simple root $\al$ and let $\Ga_{\al}$ be the stabilizer of
$\al$ in $\Ga$. Let $E_{\al}$ be the fixed field of $\Ga_{\al}$.
Since $\G$ is adjoint, for each $\al \in \Delta(\G',\T')$ we have
an associated co-root $\ca$ uniquely characterized by
$$
(\ca,\beta)={\bf 1}_{\al\beta}
$$
where $(\,,\,)$ is the natural  $\Ga$-equivariant pairing
$$
\X_*(\T')\times \X^*(\T')\ra \Z
$$
(and ${\bf 1}$ is the delta function).

Since $\al$ is defined over $K:=E_{\al}$, so is $\ca$. We have
then a morphism
$$
\ca\,:\, \gm\ra \T',
$$
defined over $K$ and consequently a continuous homomorphism
$$
\ca_{\A}\,:\, \gm(\A_K)\ra \T'(\A_K).
$$
Let $\phi_{\al}$ be the composite
$$
\phi_{\al} \,:\, \gm(\A_K) \stackrel{\ca_{\A}}{\lra}
\T'(\A_K)\stackrel{\Nm_{\A_K/\A_F}}{\lra} \T'(\A_F).
$$
By Section~\ref{sect:norms},
$$
\phi_{\al} \,:\, \gm(K)\ra \T'(F).
$$
Thus, if $\chi$ is a character of $\T'(\A_F)/\T'(F)$ then
$\xi_{\al} =\chi\circ \phi_{\al}$ is an automorphic character of
$\gm$ over $K$, i.e., a Hecke character. Write $\xi_\alpha =
\prod_{w \in \Val(K)} \xi_{\alpha, w}$.

Let $v$ be a place of $F$, $w$ a place of $K$ lying over $v$ and
$u$ a place of $E$ lying over $w$. Let $F_v\subset K_w\subset E_u$
be the corresponding completions and $\iota_w,\iota_w,\iota_u$ the
respective embeddings, with $\iota_u|_K=\iota_w,$ $\iota_w|_F
=\iota_v$. Write $\X_{E_u}^*(\T'(F_v))=\X^*(\T'(F_v))$ and let
$$
\iota^*\,:\, \X^*(\T')\stackrel{\sim}{\lra} \X^*(\T'(F_v))
$$
be the isomorphism induced by $\iota_u$. As above, we have a
restriction map
$$
r_v\,:\, \Phi(\G'(F_v),\T'(F_v))\ra \Phi(\G'(F_v),\sS_v'(F_v))
$$
mapping $\Delta(\G'(F_v),\T'(F_v))$ to
$\Delta(\G'(F_v),\sS_v'(F_v))$. The local Galois group $\Ga_{u/v}$
acts transitively on the fibers of
$$
r_v\,:\, \Delta(\G'(F_v),\T'(F_v))\ra
\Delta(\G'(F_v),\sS_v'(F_v)).
$$
Next suppose that $K_w/F_v$ is unramified. Suppose also that
$\xi_{\al,w}$ is unramified. Let
$$
L_w(s,\xi_{\al,w}) =(1-\xi_{\al,w}(\varpi_w)q_w^{-s})^{-1}
$$
where $\varpi_w$ is a prime element for $K_w$. Let
$$
\al_u=\iota_u^*(\al), \,\,\, \vartheta_v =r_v(\al_u)\in
\Delta(\G'(F_v),\sS'_v(F_v)).
$$ Let $\check{\vartheta}_v\in
\X_*(\sS'_v(F_v))$ be the associated co-character. Let
\begin{align*}
\ell(\vartheta_v) & := \# \{ \beta\in \Delta(\G'(F_v),\T'(F_v))
\,\,\, \text{ with }
\,\,\, r_u(\beta)=\vartheta_u \}\\
               & = \# \{ \Ga_{u/v}-\text{orbit of } \,\,\, \al_u \,\,\,
\text{ in } \,\,\, \Delta(\G'(F_v),\T'(F_v))\}
\end{align*}

\begin{prop}
\label{prop:41.1} With the above notations, we have
$$
L_w(s,\xi_{\al,w})=(1-\chi_v(\check{\vartheta}_v(\varpi))\cdot
q_v^{-\ell(\vartheta_v)s})^{-1}.
$$
\end{prop}

\begin{proof}
It suffices to prove
\begin{itemize}
\item $\ell(\vartheta_u) =[K_w:F_v]$; \item $\xi_{\al,w}
(\varpi_w) =\chi_v(\check{\vartheta}_u(\varpi_w))$.
\end{itemize}
let $\Ga_{u/w}=\Gal(E_u/K_w)$. Then the fixer of $\al_u$ in
$\Delta(\G'(F_v),\T'(F_v))$ is $\Ga_{u/w}$. Thus
$$
\ell(\vartheta_u) =[\Ga_{u/v} :\Ga_{u/w}] = [K_w:F_v].
$$
For the the second assertion if suffices to prove that
$$
\Nm_{K_w/F_v} (\ca_w(\la)) =\check{\vartheta}_u(\la), \,\,\,
\la\in F_v^\times.
$$
The map $\iota_u$ induces an isomorphism
$$
\iota_u^*\,:\, \X_*(\T')\stackrel{\sim}{\lra} \X_*(\T'(F_v)).
$$
Then $\iota_u^*$ preserves the natural pairing of roots and
co-characters. The co-character $\ca_w$ is obtained from $\ca$ by
base extension to $K_w$:
\begin{align*}
\ca\,:\,   &  \gm/K\ra \T'/K \\
\ca_w\,:\, &   \gm/K_w\ra \T'/K_w
\end{align*}
We have for $\la\in F_v^*$
\begin{align*}
c & = \Nm_{K_w/F_v} (\ca_w(\la)) \\
  & = \prod_{\sigma\in \Ga_{u/v}} \sigma(\ca_w(\sigma^{-1}(\la))) \\
  & = \prod_{\sigma\in \Ga_{u/v}/\Ga_{u/w}} \sigma \circ (\ca_w(\la))
\end{align*}
Here we regard $\ca_w\in \X_*(\T'(F_v))$ and $\sigma \circ$ is the
natural action of $\Ga_{u/v} $ on $\X_*(\T'(F_v))$. Since $c$ is
fixed by $\Ga_{u/v}$, we have $c\in \X_*(\sS_v(F_v))$. In fact,
the morphism $\gm\stackrel{c}{\lra} \T'(F_v)$ is defined over
$F_v$. Let $\tilde{\S}(F_v)$ be the image of $c$ in $\T'(F_v)$.
Then $\tilde{\S}(F_v)\cdot \sS_v(F_v) =\sS_v(F_v)$, since $\sS'_v$
is a maximal split $F_v$-torus in $\T'(F_v)$. Thus
$\tilde{\S}(F_v)\subset \sS_v$ and we have $\gm\stackrel{c}{\lra}
\sS'_v$.

It will suffice to prove that $c=\check{\vartheta}_u$, as elements
of $\X_*(\sS'_v(F_v))$:
$$
(c,\gamma)_{\sS'_v} =(\check{\vartheta}_u, \gamma)_{\sS'_v}
$$
for all $\gamma\in \Delta(\G'(F_v), \sS'_v(F_v))$. For this it
suffices to prove that
\begin{itemize}
\item $(c,\vartheta_u)_{\sS'_v} =1$; \item $(c,\gamma)_{\sS'_v}=
0$ for all $\gamma\in \Delta(\G'(F_v),\sS'_v(F_v))$ with $\gamma
\neq \vartheta_u$.
\end{itemize}
We have
$$
(c,\vartheta_v)_{\sS'_v} = (c, \al_u)_{\T'(F_v)} =\sum_{\sigma\in
\Ga_{u/v}/\Ga_{u/w}} (\sigma\circ \ca_w, \al_u)_{\T'(F_v)}.
$$
Next we have $\ca_w=\ca_u$ in $\X_*(\T'(F_v))$. Here
$$
\ca_w\,:\, \gm/E_u\ra \T'/E_u
$$
is obtained by base extension from
$$
\ca \,:\, \gm/E\ra \T'/E,
$$
and
$$
\al_u\,:\, \T'/E_u\ra \gm /E_u
$$
is obtained by base extension from
$$
\al\,:\, \T'/E\ra \gm/E.
$$
Since extending the scalars preserves the natural pairing between
roots and co-roots, we have the claim.
\end{proof}

\subsection{Infinite products}\label{sect:expression-infinite} The next step is to express an infinite product of the
form
\begin{equation}\label{eq:expression-infinite}
\prod_{v \in \Val(F)} \prod_{\vartheta\in
\Delta(\G'(F_v),\sS'_v(F_v))}
(1-\chi_v(\check{\vartheta}(\varpi_w))
q_v^{-\ell_v(\vartheta)s_\vartheta})^{-1}
\end{equation}
in terms of Hecke $L$-functions; this will be important in
the regularization of height integrals. The above expression, as
written, has an infinite number of complex variables.
In practice, the number of variables is finite, since
the complex numbers $s_\vartheta$ are Galois invariant
in the following sense.

Consider the commutative diagram

\

\centerline{ \xymatrix{
 \X^*(\T') \ar[r]_{\iota_v^*}  \ar[d]^{r} & \X^*(\T(F_v)) \ar[d]^{r_v}  \\
 \X^*(\sS') \ar[r]_{\iota_v^*}  & \X^*(\sS_v(F_v)).
} }

\

\noindent Let $\al\in \Delta(\G',\T')$ and let $\mO =\Ga\cdot \al$
be the orbit of $\al$ (for the action of $\Ga$ on
$\Delta(\G',\T')$). We set further
$\mO_v=r_v(\iota_u^*(\mO))\subset \Delta(\G'(F_v),\sS'_v(F_v))$.
For each $\vartheta\in \mO_v$ choose a $\beta\in \Delta(\G',\T')$
so that $r_u(\iota_u^*(\beta))=\vartheta$. We now require that
$s_\vartheta$ depend only on the Galois orbit of $\beta$ in
$\Delta(\G', \T')$, i.e., only on $\mO$; we will denote the common
value of all such $s_\vartheta$ by $s_\mO$. With the assumption,
the number of genuine complex parameters appearing in
\eqref{eq:expression-infinite} is the number of distinct Galois
orbits $\mO$ in $\Delta(\G', \T')$.

\

Fix an orbit $\mO$. For any $\beta \in \mO$, we have a field
$E_\beta$ over which $\beta$ is defined. We have already described
how to associate to $\beta$ a Hecke character
$\xi_\beta(\chi)$ of $E_\beta$. The
Hecke $L$-function $L(s, \xi_\beta(\chi))$ depends only on the
Galois orbit $\mO$, and not on the particular $\beta$. For this
reason, we denote the $L$-function $L(s, \xi_\beta(\chi))$ by
$L(s, \xi_\mO(\chi))$. An argument similar to the proof of Theorem
4.1.3 of \cite{BT2} leads to the following proposition:

\begin{prop}\label{Hecke-product}
We have
\begin{equation}
\prod_{v \in \Val(F)} \prod_{\vartheta\in
\Delta(\G'(F_v),\sS'_v(F_v))}
(1-\chi_v(\check{\vartheta}(\varpi_w))
q_v^{-\ell_v(\vartheta)s_\vartheta}) = \prod_{\mO} L(s_\mO,
\xi_\mO(\chi)),
\end{equation}
where the latter product is over all Galois orbits in $\Delta(\G',
\T')$.
\end{prop}

\section{Eisenstein series and spectral theory}
\subsection{Basic spectral theory}\label{spectral-theory}

Let $\B$ be a minimal parabolic subgroup of $\G$. We will work only
with standard parabolic subgroups and up to association. A typical
parabolic subgroup is denoted by $\P$. We denote the Levi factor of
$\P$ by $\M_\P$, its unipotent radical by $\sN_\P$, and the split
component of the center of $\M_\P$ by $\sA_\P$. We let $X(\M_\P)_\Q$
be the group of characters of $\M_\P$ defined over $\Q$, and we set
$\ma_\P = \Hom(X(\M_\P)_\Q, \R)$. Clearly, $\ma_\P^* = X(\M_\P)_\Q
\otimes \R$. We denote by $\Delta_\P$ the set of the simple roots of
$(\P, \sA)$. We usually drop the subscript $\P$.

For $m = (m_v)_v \in \M(\A)$ define a vector $H_\M(m) \in \ma_\P$
by
\begin{equation}
e^{\langle H_\M(m) , \chi \rangle} = |\chi(m)| = \prod_v
|\chi(m_v)|_v.
\end{equation}
for all $\chi \in X(\M)_\Q$. This is a homomorphism
\begin{equation}
\M(\A) \longrightarrow \ma_\P.
\end{equation}
We let $\M(\A)^1$ be the kernel. Then
\begin{equation}
\M(\A) =\M(\A)^1 \times A(\R)^0.
\end{equation}
By Iwasawa decomposition, any $x \in \G(\A)$ can be written as
$nmak$ with $n \in \sN(\A), m \in \M(\A)^1, a \in A(\R)^0, k \in
\K$. Set $H_\P(x) := H_\M(a) \in \ma_\P$. Denote the restricted Weyl
group of $(\G, \sA)$ by ${\mathcal W}$. The group ${\mathcal W}$ acts on $\ma_\B$ and
$\ma_\B^*$. For any $\fs\in {\mathcal W}$, fix a representative $w_{\fs}$ in the
intersection of $\G(\Q)$ with the normalizer of $\sA_\B$. If $\P_1,
\P_2$ are parabolic subgroups, we let ${\mathcal W}(\ma_1, \ma_2)$ be the
set of distinct isomorphisms $\ma_1 \to \ma_2$ obtained by
restricting elements of ${\mathcal W}$ to $\ma_1$. The groups $\P_1, \P_2$ are called
{\em associated} if ${\mathcal W}(\ma_1, \ma_2)$ is not empty. We usually
think of ${\mathcal W}(\ma_1, \ma_2)$ as a subset of ${\mathcal W}$. We let
$n(\sA)$ be the number of chambers in $\ma$.

Set
\begin{equation}
\ma_\P^+ = \{ H \in \ma_\P; \alpha(H) > 0 \text{ for all } \alpha
\in \Delta_\P\},
\end{equation}
and
\begin{equation}
(\ma_\P^*)^+ = \{ \Lambda \in \ma_\P^*; \Lambda(\check{\alpha}) > 0
\text{ for all } \alpha \in \Delta_\P\}.
\end{equation}
There is a vector $\rho_\P \in(\ma_\P^*)^+$ such that
\begin{equation}
\delta_\P(p) = |\det Ad\, p|_{\mn_\P(\A)}| = e^{2\rho_\P(H_\P(p))}
\end{equation}
for all $p\in \P(\A)$. We take the normalization of Haar measures to
be as in \cite{Arthur}.

Let $\sL^2_{cusp}(\M(\Q) \backslash \M(\A)^1)$ be the space of
functions $\phi$ in $\sL^2(\M(\Q) \backslash \M(\A)^1)$ such that
for any parabolic $\P_1 \nsubseteq \P$ we have
\begin{equation}
\int_{\sN_1(\Q) \cap \M(\Q) \backslash \sN_1(\A) \cap \M(\A)}
\phi(nm) \, dn = 0
\end{equation}
for almost all $m$. It is known that
\begin{equation}
\sL^2_{cusp}(\M(\Q) \backslash \M(\A)^1) = \bigoplus_\rho V_\rho
\end{equation}
where $V_\rho$ is the $\rho$-isotypic component of $\rho$ consisting
of finitely many copies of $\rho$ (possibly zero). The pairs $(\M,
\rho)$ and $(\M', \rho')$ are considered equivalent if there is $\fs
\in {\mathcal W}(\ma, \ma')$ such that the representation
\begin{equation}
(\fs\rho)(m') = \rho(w_{\fs}^{-1} m' w_{\fs}) \,\,\,\,\,(m' \in \M'(\A)^1)
\end{equation}
is unitarily equivalent to $\rho'$. Let $\mX$ be the set of
equivalence classes of such pairs. For any $X \in \mX$ we have a
class $\mathcal{P}_X$ of associated parabolic subgroups. If $\P$ is
any parabolic and $X \in \mX$, set
\begin{equation}
\sL^2_{cusp}(\M(\Q) \backslash \M(\A)^1)_X = \bigoplus_{(\rho: (\M,
\rho) \in X)}V_\rho.
\end{equation}
This is a closed subspace of $\sL^2_{cusp}(\M(\Q) \backslash
\M(\A)^1)$, and empty if $\P \notin \mathcal{P}_X$.

Fix $\P$ and $X \in \mX$. Suppose there is a $\P_1 \in
\mathcal{P}_X$ such that $\P_1 \subset \P$. Let $\psi$ be a smooth
function on $\sN_1(\A)\M_1(\Q) \backslash \G(\A)$ such that
\begin{equation}
\Psi_a(m, k) = \psi(amk) \,\,\,\,\,(k \in \K, m \in
\M_1(\Q)\backslash \M_1(\A), a \in \sA_1(\Q)\backslash \sA_1(\A))
\end{equation}
vanishes outside a compact subset of $\sA_1(\Q)\backslash
\sA_1(\A)$, transforms under $\K_\R$ according to an irreducible
representation, and as a function of $m$ belongs to
$\sL^2_{cusp}(\M(\Q) \backslash \M(\A)^1)$. Then the function
\begin{equation}
\hat{\psi}^\M(m) = \sum_{\delta\in \P_1(\Q) \cap \M(\Q) \backslash
\M(\Q)} \psi(\delta m)\,\,\,\,\,\, (m \in \M(\Q) \backslash
\M(\A)^1)
\end{equation}
is square-integrable on $\M(\Q) \backslash \M(\A)^1$. Define
$\sL^2(\M(\Q) \backslash \M(\A)^1)_X$ to be the span of all such
$\hat{\psi}^\M$. If no such $\P_1$ exist, latter space is the zero space.

By a result of Langlands,
\begin{equation}
\sL^2(\M(\Q) \backslash \M(\A)^1) = \bigoplus_X \sL^2(\M(\Q)
\backslash \M(\A)^1)_X.
\end{equation}
For any $\P$, let $\Pi(\M)$ denote the set of equivalence classes of
irreducible unitary representations of $\M(\A)$. For $\zeta \in
\ma_\C^*$ and $\pi \in \Pi(\M)$ let $\pi_\zeta$ be the product of
$\pi$ with the quasi-character
\begin{equation}
x \mapsto e^{\zeta H_\P(x)} \,\,\,\,\,\,(x \in \G(\A)).
\end{equation}
If $\zeta \in i\ma^*$, $\pi_\zeta$ is again unitary. This means that
$\Pi(\M)$ is a differentiable manifold which carries an action of
$i\ma^*$. We use this action to define a measure $d\pi$ on
$\Pi(\M)$.

For $\pi \in \Pi(\M)$ we let $\mathcal{H}^0_\P(\pi)$ be the space of
smooth functions
\begin{equation}
\phi: \sN(\A) \M(\Q) \backslash \G(\A) \to \C
\end{equation}
satisfying
\begin{enumerate}
\item $\phi$ is right $\K$-finite;
\item for every $x \in \G(\A)$ the function $$m \mapsto
\phi(mx)\,\,\,\,(m \in \M(\A))$$ is a matrix coefficient of $\pi$;
\item $\| \phi \|^2 = \int_\K \int_{\M(\Q)\backslash \M(\A)^1}
|\phi(mk)|^2 \, dm \, dk < +\infty$.
\end{enumerate}
Let $\mathcal{H}_\P(\pi)$ be the completion. If $\phi \in
\cH_\P(\pi)$ and $\zeta \in \ma_\C^*$ set
\begin{equation}
\phi_\zeta(x)  = \phi(x) e^{\zeta(H_\P(x))} \,\,\,\,\,(x \in \G(\A))
\end{equation}
and
\begin{equation}
(I_\P(\pi_\zeta, y)\phi_\zeta)(x) = \phi_\zeta(xy)
\delta_\P(xy)^{\frac{1}{2}} \delta_\P(x)^{-\frac{1}{2}}.
\end{equation}
Then $I_\P(\pi_\zeta)$ is a unitary representation if $\zeta \in i\ma^*$.

Given $X \in \mX$, let $\cH_\P(\pi)_X$ be the closed subspace of
$\cH_\P(\pi)$ consisting of those $\phi$ such that for all $x$ the
function $m \mapsto \phi(mx)$ belongs to $\sL^2(\M(\Q) \backslash
\M(\A)^1)_X$. Then
\begin{equation}
\cH_\P(\pi) = \bigoplus_{X} \cH_\P(\pi)_X.
\end{equation}
Let $\K_0$ be an open-compact subgroup of $\G(\A_f)$ and $W$ an
equivalence class of irreducible representations of $\K_\R$. Let
$\cH_\P(\pi)_{X, K_0}$ be the subspace of functions in
$\cH_\P(\pi)_X$ which are invariant under $\K_0 \cap \K$. Also let
$\cH_\P(\pi)_{X, K_0, W}$ be the space of those functions in $\cH_\P(\pi)_{X,
K_0}$ which transform under $\K_\R$ according to $W$. It is a
theorem of Langlands that each of the spaces $\cH_\P(\pi)_{X, K_0,
W}$ is finite-dimensional. We fix an orthonormal basis
$\cB(\pi)_X$ for $\cH_\P(\pi)_X$, for each $\pi$ and each $X$,  such
that for all $\zeta \in i \ma^*$ we have
\begin{equation}
\cB(\pi_\zeta)_X = \{ \phi_\zeta : \phi \in \cB_\P(\pi)_X \}
\end{equation}
and such that every $\phi \in\cB_\P(\pi)_X$ belongs to one of the
spaces $\cH_\P(\pi)_{X, K_0, W}$.

Suppose that $\pi \in \Pi(\M)$, $\phi \in \cH_\P^0(\pi)$, $\zeta \in
\ma_\C^*$. For $\Re(\zeta) \in \rho_\P + (\ma^*)^+$ we set
\begin{equation}
E(x, \phi, \zeta) = \sum_{\delta \in \P(\Q) \backslash \G(\Q)}
\phi_\zeta(\delta x) \delta_\P(\delta x)^{\frac{1}{2}}.
\end{equation}
For $\fs \in {\mathcal W}(\ma, \ma')$, we define the global intertwining
operator
\begin{equation}\begin{split}
(M(\fs, \pi, \zeta)& \phi_\zeta)(x) \\
& = \int_{\sN'(\A) \cap w_{\fs} \sN(\A) w_{\fs}^{-1} \backslash \sN'(\A)}
\phi_\zeta(w_{\fs}^{-1} n x) \delta_\P(w_{\fs}^{-1} n x)^{\frac{1}{2}}
\delta_\P(x)^{-\frac{1}{2}}\, dn.
\end{split}\end{equation}
Both $E(x, \phi, \zeta)$ and $M({\fs}, \pi_\zeta)\phi_\zeta$ can be analytically
continued to meromorphic functions in $\zeta$ to $\ma_\C^*$. For
$\zeta \in i\ma^*$, $E(x, \phi, \zeta)$ is a smooth function of $x$
and $M(\fs, \pi_\zeta)$ is unitary from $\cH_\P(\pi_\zeta)$ to
$\cH_{\P'}(\fs\pi_\zeta)$.


\subsection{Spectral expansion}
Let $f$ be a function on $\G(\Q) \backslash \G(\A)$. By
the {\em spectral expansion} of $f$ we mean
\begin{equation}\begin{split}\label{spectral-identity-general}
S(f, x) = \sum_{X \in \mX}\sum_{\P} & n(\sA)^{-1} \\
& \int_{\Pi(\M)}\sum_{\phi \in \cB_\P(\pi)_X}  E(x, \phi)\left(
\int_{\G(\Q) \backslash \G(\A)} \overline{E(y, \phi)}f(y) \, dy
\right)\, d\pi.
\end{split}\end{equation}
For $\M = \G$, the integrals are interpreted appropriately to give a
discrete sum. It is an
interesting problem to determine under which conditions on the
function $f$, we have $f(x) = S(f, x)$. Also, for the
application we have in mind, we need to know that the identity $f =
S_f$ holds as an identity not just of $\sL^2$-functions, but of
continuous functions.
In order to show that for a given $f$,
the spectral identity holds, we just need to prove that that right
hand side of the \eqref{spectral-identity-general} is uniformly
convergent on compact sets. This would suffice since for
all pseudo-Eisenstein series $\theta_\phi$ (II.1.10 of
\cite{MW}) we have
\begin{equation}
(S(f), \theta_\phi) = (f, S(\theta_\phi)),
\end{equation}
and since pseudo-Eisenstein series are rapidly decreasing
(Proposition II.1.10 \cite{MW}), we have $\theta_\phi =
S(\theta_\phi)$ by the spectral decomposition \cite{Arthur}.
Consequently $(f- S(f), \theta_\phi)=0$ for all pseudo-Eisenstein
series $\theta_\phi$. Then by the density theorem Theorem II.1.12 of
\cite{MW} we get the required identity. In our applications, we
use a slightly different re-arrangement of the terms of $S(f)$ as
suggested by the definition of $I(S, f, x, y)$ on page 930 of
\cite{Arthur}.

Let $f$ be a smooth function on $\G(\A)$, and suppose that $f$ is
right invariant under a compact-open subgroup $\K$ of $\G(\A_f)$.
Define a function on $\G(\Q) \backslash \G(\A)$ by
$$
F(g) := \sum_{\gamma \in \G(\Q)} f(\gamma g).
$$
Suppose that $f$ is such that the function $F$ is convergent for
all $g$, is smooth and bounded. By Lemma 4.1 of \cite{Arthur},
given $m$, there is an $n$ and functions $f_1 \in
C_c^m(\G(\R))^{K_\R}$ and $f_2 \in C_c^\infty(\G(\R))^{K_\R}$ such
that $\Delta^n*f_1 + f_2$ is the delta distribution at the
identity of $\G(\R)$. Here $\Delta$ is an appropriate element of
the enveloping algebra chosen as in the proof of Lemma 4.1 of
\cite{Arthur}. This implies that after taking the finite places
into account, we can find functions $\eta_1, \eta_2$ of compact
support on $\G(\A)$, with appropriate differentiability at the
archimedean place, such that
\begin{equation}\label{equation-reference-for-arthur}
R(\Delta^n * \eta_1 + \eta_2) F = F,
\end{equation}
where $R$ is the right regular convolution action. We get
\begin{equation}
R(\eta_1) \Delta^n F + R(\eta_2) F = F.
\end{equation}
Consider the following modified spectral expansion
\begin{equation}\begin{split}\label{spectral-identity-modified}
S'(F, x) = \sum_{X \in \mX}\sum_{\P} & n(\sA)^{-1} \\
& \int_{\Pi(\M)}\int_{\G(\Q) \backslash \G(\A)}\left(\sum_{\phi \in
\cB_\P(\pi)_X}  E(x, \phi) \overline{E(y, \phi)}\right)F(y) \, dy\,
d\pi.
\end{split}\end{equation}
We get
\begin{align*}
S'(F, x) = \sum_{X \in \mX}\sum_{\P} & n(\sA)^{-1} \\
& \int_{\Pi(\M)}\int_{\G(\Q) \backslash \G(\A)}\left(\sum_{\phi \in
\cB_\P(\pi)_X}  E(x, \phi) \overline{E(y, \phi)}\right)F(y) \, dy\,
d\pi \\
= \sum_{X \in \mX}\sum_{\P} & n(\sA)^{-1} \int_{\Pi(\M)}\int_{\G(\Q)
\backslash \G(\A)}\left(\sum_{\phi \in
\cB_\P(\pi)_X}  E(x, \phi) \overline{E(y, \phi)}\right)\\
& (R(\eta_1) \Delta^n F (y)+ R(\eta_2) F(y) )\, dy\, d\pi.
\end{align*}
Consequently, $S'(F, x)$ will be the sum of the following two series
\begin{align*}
S_1(F, x) = \sum_{X \in \mX}\sum_{\P} & n(\sA)^{-1}
\int_{\Pi(\M)}\int_{\G(\Q) \backslash \G(\A)}\left(\sum_{\phi \in
\cB_\P(\pi)_X}  E(x, \phi) \overline{E(y, I_\P(\pi, \tilde{\eta}_1)\phi)}\right)\\
& \Delta^n F (y)\, dy\, d\pi
\end{align*}
and
\begin{align*}
S_2(F, x) = \sum_{X \in \mX}\sum_{\P} & n(\sA)^{-1}
\int_{\Pi(\M)}\int_{\G(\Q) \backslash \G(\A)}\left(\sum_{\phi \in
\cB_\P(\pi)_X}  E(x, \phi) \overline{E(y, I_\P(\pi, \tilde{\eta}_2)\phi)}\right)\\
& F (y)\, dy\, d\pi
\end{align*}
for appropriately chosen compactly supported function
$\tilde{\eta}_1, \tilde{\eta}_2$. By Lemma 4.4 of
\cite{Arthur}, there are $N, r_0$ and a continuous
seminorm $\| \cdot \|_{r_0}$ on $C_c^{r_0}(\G(\A))$ such that if $r \geq
r_0$ and $\eta$ is $\K$-finite in $C_c^r(\G(\A))$,
\begin{align*}
\sum_{X \in \mX}\sum_{\P} & n(\sA)^{-1} \int_{\Pi(\M)} \left|\sum_{\phi
\in \cB_\P(\pi)_X}  E(x, \phi) \overline{E(y,
I_\P(\pi,\eta)\phi)}\right|\, d\pi
\end{align*}
is bounded by $ \| \eta\|_{r_0}\cdot \| x \|^N\cdot \| y \|^N$. Here we have
defined the height $\|\cdot\|$ on $\G(\A)$ as in \cite{Arthur}, page
918. Since $m$ can be taken to be
arbitrarily large, all we need to verify in order to get the uniform
convergence is the convergence of the following integrals
\begin{equation}\label{integral1-spectral}
\int_{\G(\A)} |f(y)|\cdot \| y\|^N \, dy
\end{equation}
and
\begin{equation}\label{integral2-spectral}
\int_{\G(\A)} |\Delta^n f(y)|\cdot \| y\|^N \, dy.
\end{equation}

We summarize this discussion:

\begin{lem}\label{lem:spectral-expansion-general}
Let $n, N$ be as above. Assuming convergence of the integrals
\eqref{integral1-spectral} and \eqref{integral2-spectral} we have
\begin{equation}
F(x) = S_1(F, x) + S_2(F, x).
\end{equation}
\end{lem}


\subsection{Truncations} Let $T \in \ma_\B^+$. Let $\Delta_0(\B, \sA)$ be a
set of simple roots of $(\B, \sA)$, and $\mathsf{d}_\B(T) :=
\min_{\alpha \in \Delta_0(\B, \sA)} \{\alpha(T)\}$. We will need
those $T$ which are \emph{sufficiently regular}; this means that
$\mathsf{d}_\B(T)$ is large. Recall that for sufficiently regular
$T$, the Arthur truncation $\wedge^T$ acts on functions on $\G(\Q)
\backslash \G(\A)$ \cite{Arthur3}.

\begin{lem}
\begin{enumerate}
\item If $\phi$ is locally bounded, then $\wedge^T \phi$ is defined
everywhere.
\item If $\phi: \G(\Q) \backslash \G(\A) \to \C$ is locally $\sL^1$,
then we have
$$
\wedge^T \wedge^T \phi(g) = \wedge^T \phi(g)
$$
for almost all $g$. If $\phi$ is locally bounded, then the identity
holds for all $g$.
\item Suppose $\phi_1, \phi_2:\G(\Q) \backslash \G(\A) \to \C$ are
locally $\sL^1$. If $\phi_1$ is of moderate growth and $\phi_2$ is
rapidly decreasing, then
$$
\int_{\G(\Q) \backslash \G(\A)} \overline{\wedge^T \phi_1(g)}
\phi_2(g) \, dg = \int_{\G(\Q) \backslash \G(\A)}
\overline{\phi_1(g)} \wedge^T \phi_2(g)\, dg.
$$
\end{enumerate}
\end{lem}
Observe that the lemma implies that if $\phi_1, \phi_2$ are as
above, then
\begin{equation}
\int_{\G(\Q) \backslash \G(\A)} \overline{\wedge^T \phi_1(g)}
\phi_2(g) \, dg = \int_{\G(\Q) \backslash \G(\A)} \overline{\wedge^T
\phi_1(g)} \wedge^T \phi_2(g)\, dg.
\end{equation}

Given $\P$, and $\pi \in \Pi(\M_{\P}(\sA))$ and $\lambda \in i \ma_{\P}^*$
we define an operator
\begin{equation}
\Omega_{X, \pi}^T(\P, \lambda) : \mathcal{H}_\P^0(\pi)_X \to
\mathcal{H}_\P^0(\pi)_X
\end{equation}
by setting
\begin{equation}
(\Omega_{X, \pi}^T(\P, \lambda)\pi, \pi') = \int_{\G(\Q) \backslash
\G(\A)^1} \wedge^T E(x, \phi, \lambda)\overline{\wedge^T E(x, \phi',
\lambda)}\, dx,
\end{equation}
for each pair of vectors $\phi, \phi' \in \mathcal{H}_\pi^0(\pi)_X$.

Next let $\K_0$ be a subgroup of finite index in $\K_f$. Suppose
that $W$ is a finite dimensional representation of $\K_\infty$.
Given $\P \supset \B$, $\pi \in \P(\M_\P(\A)^1)$, let
$\mathcal{H}_\P^0(\pi)_X^{\K_0}$  be the space of $\K_0$-invariant
functions in $\mathcal{H}_\P^0(\pi)_X$, and let
$\mathcal{H}_\P^0(\pi)_X^{\K_0, W}$ be the subspace of functions
in $ \mathcal{H}_\P^0(\pi)_X^{\K_0}$ which transform under
$\K_\infty$ according to $W$. If a linear operator $A$ on
$\mathcal{H}_\P^0(\pi)_X$ leaves any of these subspaces
$\mathcal{H}_\P^0(\pi)_X^\square$, with $\square = \K_0$ or
$\square= \K_0, W$, invariant, we denote by $A_*$ the restriction
of $A$ to the appropriate subspace. Choose the differential
operator $\Delta$ as in the proof of Lemma 4.1 of \cite{Arthur}.
Then $\Delta$ acts on $\mathcal{H}_\P^0(\pi)_X$ through each of
the representations $I_\P(\pi_\zeta)$, and we denote the action by
$I_\P(\pi_\zeta, \Delta)$. This action leaves the two subspaces
mentioned above invariant. For $\zeta \in i\ma_{\P}^*$, the action
of $I_\P(\pi_\zeta, \Delta)$ is via a scalar greater than 1.

\begin{lem}[Arthur \cite{Arthur2}]\label{arthur-main}
There exist integers $C_0, d_0$,
and $m$ such that for any subgroup $\K_0
\subset \K_f$ of finite index, and any $T \in i\ma_\B^+$ with
$\mathsf{d}_{\B}(T) > C_0$, the expression
\begin{equation}
\sum_{X \in \mathfrak{X}} \sum_\P \sum_{\pi \in \Pi(\M_\P(\A)^1)}
|n(\sA)|^{-1} \int_{i \ma_\P^*/ i \ma_\G^*} \| \Omega_{X, \pi}^T(\P,
\zeta)_{\K_0}. I_\P(\pi_\zeta, \Delta^m)^{-1}\|_1 d \zeta
\end{equation}
is bounded by $C_{\K_0} (1 + \| T\|)^{d_0}$. Here $C_{\K_0}$ is a
constant which depends only on $\K_0$, and $\|\cdot\|_1$ is the trace
class norm.
\end{lem}

\subsection{A bound for Eisenstein series} We recall a bound
for Eisenstein series embedded in the proof of Proposition 2 of
\cite{Lapid}. Fix a large compact set
$C$ of $\G(\A)$ and a small compact open subgroup $\K$ of $\G(\A_f)$,
and let $f$ be a sufficiently differentiable function with support
in $C$.
\begin{lem} Fix a compact set $C'$ in $\G(\A)$. There are constants $c, N>0$
depending only on the support of $f$ such that for sufficiently
regular $T$ we have
\begin{equation}
|E(g, I(f, \pi_\lambda)\varphi, \lambda)| \leq c \| g \|^N
\|f\|_\infty\cdot \|\wedge^T E(\cdot, \varphi, \lambda) \|_{\L^2(\G(\Q)
\backslash \G(\A)) }.
\end{equation}
\end{lem}
This lemma combined with Lemma \ref{arthur-main} together with Lemma
4.1 of \cite{Arthur}
implies the following result.

\begin{prop}\label{important-proposition}
For $\phi \in \cB_\P(\pi)_X$ let $\Lambda(\phi)$
be defined by $\Delta\cdot \phi = \Lambda(\phi)\cdot \phi$. Then there is
an $l>0$ such that
\begin{equation}
\sum_{X \in \mX}\sum_{\P} n(\sA)^{-1} \int_{\Pi(\M)}\left( \sum_{\phi
\in \cB_\P(\pi)_X} \Lambda(\phi)^{-m}|E(e, \phi)|^2 \right)\, d\pi
\end{equation}
is convergent. The outermost summation is only over those classes
for which $\Lambda(\phi) \ne 0$ are fixed by $\K_0$.
\end{prop}

\section{Spherical functions and bounds for matrix coefficients}
\subsection{Spherical functions}\label{spherical-functions}
Let $\pi$ be an infinite-dimensional automorphic representation of
$\G$, and suppose $\phi$ is a right $\K$-finite automorphic form in
the space of $\pi$. Here $\K = \prod_v \K_v$. Define a function
$M(\phi, g)$ on $\G(\A)$ by
\begin{equation*}
M(\phi, g) = {1 \over (\text{vol }\sK)^2} \int_{\sK}\int_{\sK}
\phi(\kappa g \kappa') \, d\kappa \, d\kappa'.
\end{equation*}
If $\pi$ has no $\K$-invariant vectors, then the above integral is
zero, and we may assume that $\phi$ is right $\sK$ invariant.
Let $\V_\pi$ be a vector space on which $\G(\A)$ acts via a
representation which is isomorphic to the space of $\K$-finite
vectors. Let $j\,:\, \V_\pi \rightarrow \sH_\pi$ be the
intertwining map. We have used the same notation for a
representation and its underlying representations space. By a
standard result of Jacquet, Langlands, and Flath, we have
\begin{equation*}
\V_\pi = \bigotimes_v \V_{\pi,v}.
\end{equation*}
Since $\phi$ is right $\sK$-finite, it is in the image of $j$. Let
$w= j^{-1}(\phi)$. Then

\begin{equation*}
w \in \V_\pi^{\sK} = \left( \bigotimes_{v \notin S_\infty} \V_{\pi,
v}^{\sK_v}\right) \otimes \left( \bigotimes_{v \in S_\infty}
\V_{\pi,v} \right).
\end{equation*}
Set $S = S_F$. We know that for $v \notin S$ we have ${\rm dim}\, \V_{\pi,v}^{\sK_v} =1$. Fix a non-zero element $e_v$ in each of
these spaces. This means then
\begin{equation*}
w = (\otimes_{v \notin S} e_v) \otimes w_S,
\end{equation*}
where
$$
w_S \in \left(\bigotimes_{v \in S \backslash S_\infty}
\V_{\pi,v}^{\sK_v}\right) \otimes \left( \bigotimes_{v \in S_\infty}
\V_{\pi,v}\right).
$$
Next we examine $M(\phi, g)$. Define a functional $\lambda$ on
$\V_\pi$ by
\begin{equation*}
\lambda (\nu) = \int_{\sK} j(\nu) (\kappa) \, d\kappa,
\end{equation*}
(for $\nu \in \V_\pi$). Then it is easily seen that
$$
\lambda \in \bigotimes_{v \notin S } \widetilde{\V}_{\pi,v}^{\sK_v}
\otimes \V_S^*.
$$
Here $\V_S = \bigotimes_{v \in S} \V_{\pi,v}$, $\V_S^*$ is the dual
space and $\widetilde{\V}_\pi$ is the smooth dual of the local
representation $\V_\pi$. Since the smooth dual of an admissible
representation is admissible, it follows that
$\dim\widetilde{\V}_{\pi,v}^{\sK_v} =1$. For each $v \notin S$
choose an element $\xi_v$ in this space in such a way that
$\xi_v(e_v) =1$. This then means that, as in the case of $w =
j^{-1}(\phi)$, we have
\begin{equation*}
\lambda= ( \otimes_{v \notin S} \xi_v) \otimes \lambda_S,
\end{equation*}
with $\lambda_S$ in the obvious space. Combining this identity with
the similar identity for $w$ we obtain
\begin{equation*}\begin{split}
\lambda(v) & = (\prod_{v \notin S} \xi_v(e_v))\cdot \lambda_S(w_S) \\
& = \lambda_S(w_S).
\end{split}\end{equation*}
This implies that
\begin{equation*}\begin{split}
\lambda_S (w_S) & = \int_{\sK} j(w)(\kappa) \, d\kappa \\
& = \int_{\sK^S} j(w)(\kappa) \, d\kappa.
\end{split}\end{equation*}
Here $\sK^S = \prod_{v \in S} \sK_v$. Consider an embedding $\eta :
\G(\A_S) \longrightarrow \G(\A)$ given by
\begin{equation*}
g \mapsto (1, 1, \dots, 1, g).
\end{equation*}
Obviously,
\begin{align*}
{\rm vol}(\sK^S)M(\phi, g) & = \lambda(\pi (g) w)  \\
   & = \prod_{v \notin S}\varphi_{v, \pi}(g_v) \cdot\int_{\sK^S}
\phi( \kappa\eta(g_S)) \, d\kappa,
\end{align*}
where $\varphi_{v, \pi}(g_v) = \xi_v(\pi_v(g_v) e_v)$ is the
normalized local spherical function.

\begin{coro}\label{coro-spherical-satake} If $\K = \prod_v \K_v$ is such that
for each local place $v$, including the archimedean places, the
local compact subgroups $\K_v$ satisfy the hypotheses I and II of
\cite{satake1}, then
\begin{equation}
M(\phi, g) = \phi(e) \prod_v \varphi_{v, \pi}(g_v)
\end{equation}
if each $\pi$ has a $\K_v$-fixed vector; otherwise it is zero.
\end{coro}
We note that local $\K_v$ satisfying I and II of
\cite{satake1} exist by \cite{GV}, for archimedean places, and
\cite{B-T, silberger, Tits}, for non-archimedean places.

\subsection{Bounds on Matrix Coefficients}
\label{sect:bounds}

In this section we recall an important result
of H. Oh that is used in estimates leading
to the proof of the spectral expansion \eqref{spectral-identity-general}.

Let $k$ be a non-archimedean local field of ${\rm char}(k) \ne 2$, and residual
degree $q$. Let $\H$ be the group of $k$-rational points of a
connected reductive split or quasi-split group with $\H / \sZ(\H)$
almost $k$-simple. Let $\sS$ be a maximal $k$-split torus, $\B$ a
minimal parabolic subgroup of $\H$ containing $\sS$ and $\sK$ a good
maximal compact subgroup of $\H$ with Cartan decomposition $\G = \K
\sS(k)^+ \K$. Let $\Phi$ be the set of non-multipliable roots of the
relative root system $\Phi(\H, \sS)$, and $\Phi^+$ the set of
positive roots in $\Phi$. A subset $\mathcal{S}$ of $\Phi^+$ is
called a strongly orthogonal system of $\Phi$ if any two distinct
elements $\alpha$ and $\alpha'$ of $\mathcal{S}$ are strongly
orthogonal, that is, neither of $\alpha \pm \alpha'$ belongs to
$\Phi$. Define a bi-$\K$-invariant function $\xi_{\mathcal{S}}$ on
$\H$ as follows: first set

\begin{equation*}
n_{\mathcal{S}}(g) = {1 \over 2} \sum_{\alpha \in \mathcal{S}}
\log_q \vert \alpha (g) \vert, \end{equation*} then
\begin{equation*}
\xi_{\mathcal{S}}(g) = q^{-n_{\mathcal{S}}(g) } \prod_{ \alpha \in
\mathcal{S}} ( \frac{ (\log_q \vert \alpha(g) \vert ) (q-1) + (q+1)
} { q+1} ).
\end{equation*}

\begin{thm}[\cite{Oh}, Theorem 1.1]
\label{Oh} Assume that the semi-simple $k$-rank of $\H$ is at least $2$ and
let $\mathcal{S}$ be a strongly orthogonal system of $\Phi$.
Then for any unitary representation $\varrho$ of $\H$ without an
invariant vector and with $\sK$-finite unit vectors $\nu$ and
$\nu'$, one has

\begin{equation*}
\vert (\varrho(g)\nu,\nu') \vert \leq (\dim(\sK\nu)
\dim(\sK\nu'))^{1 \over 2}\cdot \xi_{\mathcal{S}}(g),
\end{equation*}
for any $g \in \H$.
\end{thm}

\



\begin{coro}
\label{n} Let $v \notin S_F$. Let $\varphi_v$ be the normalized
spherical function associated with an infinite dimensional
unramified principal series representation of $\G(F_v)$. Suppose
that the semi-simple rank of $\G(F_v)$ is at least 2. Then for
each $\vartheta \in \De(\G(F_v), \S_v(F_v))$ we have
\begin{equation*}
\vert \varphi_v (\check\vartheta(\varpi_v)) \vert < 2
q_v^{-{\ell_v(\vartheta)\over 2}}.
\end{equation*}
\end{coro}

\begin{proof}
Fix $\vartheta$. By definition, a singleton is a strongly orthogonal
set. Set
\begin{equation*}
\mathcal{S}_v = \{ \vartheta\}.
\end{equation*}
Then
\begin{equation*}
n_{\mathcal{S}_v}(\check\vartheta(\varpi_v)) = {\ell_v(\vartheta)
\over 2},
\end{equation*}
and
\begin{equation*}
\xi_{\mathcal{S}_p}(\check\vartheta(\varpi_v)) =
q_v^{-{\ell_v(\vartheta) \over 2}} {2q_v \over q_v+1}.
\end{equation*}
Since $\nu$ and $\nu'$ that define the spherical function are
$\sK_v$-invariant, this equation combined with Theorem \ref{Oh}
implies
\begin{equation*}\begin{split}
\vert \varphi_v (\check\vartheta(\varpi_v)) & \leq
\xi_{\mathcal{S}_v}(\check\vartheta(\varpi_v))
) \\
&\leq q^{-{\ell_v(\vartheta) \over 2}} {2q_v \over q_v+1} \\
&< 2 q_v^{-{\ell_v(\vartheta) \over 2}},
\end{split}\end{equation*}
which gives the claim.
\end{proof}

We also need a similar bound on spherical functions when the
semi-simple rank is equal to one. In this case, local considerations
do not suffice, as the trivial representation may not be isolated in
the unitary dual of the local group. However, for
our purposes it will suffice to obtain a bound
for a restricted class of representations.
Let $\pi$ be an infinite-dimensional unitary
irreducible automorphic representation of $\G$. Suppose
$\pi = \otimes_{v} \pi_v$. Extending $S_F$ if necessary, we may
assume that for $v \notin S_F$, $\pi_v$ is an unramified
representation.

\begin{prop}
\label{one-dimensional}
For all $v \notin S_F$, the representation
$\pi_v$ is infinite-dimensional.
\end{prop}

\begin{proof} Realize $\pi$ on a Hilbert subspace $\V$
of $\sL^2(\G(F) \backslash \G(\A))$.
Denote by $\V^\infty$ the subspace of $\V$ consisting of all vectors
$\nu$ such that:
\begin{itemize}
\item $\nu$ is $\K_f$-finite, and \item for all archimedean places
$v$, the map $\G(F_v) \to \V$ given by $g \mapsto \pi(g)\nu$ is
$C^\infty$.
\end{itemize}
It is standard that $\V^\infty$ is $\G(\A)$-stable, and dense in
$\V$ (with respect to the $\sL^2$-topology). The image of
$\V^\infty$ in $\sL^2$ consists of smooth functions. Suppose that
for some $v \notin S_F$ the representation $\pi_v$ is not
infinite-dimensional. Since this representation is admissible, it
must be one-dimensional. Concretely, if we pick an element $\varphi
\in \V^\infty$, we have
\begin{equation*}
\varphi(\gamma x_v) = \chi_v(x_v) \varphi(e), \quad{\gamma \in
\G(F), x_v \in \G(F_v).}
\end{equation*}
Here $\chi_v: \G(F_v) \to \C^\times$ is a one-dimensional
representation of $\G(F_v)$. Next, by Satz 6.1. of \cite{Kneser2} we
know that the commutator subgroup $\G'(\A)$ of $\G(\A)$ is contained
in the closure of $\G(F)\G(F_v)$. As $\G'(\A)$ has no non-trivial
one-dimensional representations, we conclude that for all $\varphi
\in \V^\infty$ we have
\begin{equation}\label{2}
\varphi(g') = \varphi(e),
\end{equation}
for all $g' \in \G'(\A)$. Since the subspace $\V^\infty$ is
$\G(\A)$-invariant, \eqref{2} must hold for all right translates
of $\varphi$ by elements of $\G(\A)$. Hence
\begin{equation*}
\varphi(xg) = \varphi(g), \quad{x \in \G'(\A), g \in \G(\A)}.
\end{equation*}
Since $\G'(\A)$ is a normal subgroup of $\G(\A)$, it follows that we
must also have
\begin{equation*}
\varphi(gx) = \varphi(g), \quad{x \in \G'(\A), g \in \G(\A)},
\end{equation*}
i.e., every element of $\V^\infty$ is invariant under the
restriction of the representation $\pi$ to the subgroup $\G'(\A)$.
This implies that the representation $\pi$ on $\V^\infty$ factors
through $\G(\A) / \G'(\A)$, which is an abelian group. Since
irreducible admissible representations of abelian groups are all
one-dimensional, the space $\V^\infty$ must be one-dimensional.
Finally, use the fact that $\V^\infty$ is dense in $\V$ to conclude
that $\V$ is also one-dimensional.
\end{proof}

\begin{thm}
\label{bound} There is an absolute constant $c>0$ with the
following property. Let $\pi = \otimes_v \pi_v$ be an
infinite-dimensional irreducible automorphic representation of $\G$;
for $\G= \PGL_2$ assume that $\pi$ is not the
automorphic representation associated with an Eisenstein series
induced from a Borel subgroup. Let $v$ be a place with $v\notin
S_F$, and $\varphi_v$ the normalized spherical function of $\pi_v$.
Then for all $\alpha \in \De(\G(F_v), \S_v(F_v))$ we have
\begin{equation*}
\vert \varphi_v (\check\vartheta(\varpi_v)) \vert < q^{-c
{\ell_v(\vartheta)}}.
\end{equation*}
\end{thm}

\begin{proof}
According to the result of Corollary~\ref{n},
we recognize two cases: \\

{\em The case where s.s. rank is equal to one. } In this case, $\G$
is a form of $\PGL_2$ or ${\rm PGU}_3$. For the case of $\PGL_2$, if
$\G$ is not split, then it has to be a quaternion algebra, and in
this case by the Jacquet-Langlands correspondence, there is an
irreducible cuspidal automorphic representation $\pi'$ of $\GL_2$,
such that for all $v \notin S$, we have $\pi_v = \pi'_v$ as
representations of $\GL_2(F_v)$. Since the local representation
obtained this way is unramified, there must exist a pair of
unramified quasi-characters $\chi$ and $\chi'$ of $F_v^\times$ such
that $\pi'_v= \pi(\chi, \chi')$. Since this representation has
trivial central character, we must have $\chi' = \chi^{-1}$. In this
case, if $\alpha$ is the unique positive root, we have
\begin{equation*}
t_v(\alpha) =  \begin{pmatrix} 1 \\ & \varpi_v \end{pmatrix}.
\end{equation*}
Next if $\chi\ne 1$, by the formula of Casselman
(\cite{Casselman2}), we have
\begin{equation*}\begin{split}
\varphi_v \begin{pmatrix} 1 \\ & \varpi_v \end{pmatrix} & = {q^{-{1
\over 2}} \over 1 + q^{-1}} \{ {1 - q^{-1} \chi(\varpi_v)^{-2} \over
1 - \chi(\varpi_v)^{-2}} \chi(\varpi_v) + {1 - q^{-1}
\chi(\varpi_v)^{2}
\over 1 - \chi(\varpi_v)^{2}} \chi(\varpi_v)^{-1} \} \\
&= {q^{-{1 \over 2}} \over 1 + q^{-1}} (\chi(\varpi_v) +
\chi(\varpi_v)^{-1}).
\end{split}\end{equation*}
The same formula holds for $\chi =1$ by analytic continuation. Since
$\chi$ is an unramified character, we have $\chi = \vert \cdot\vert
^s$, for some complex number $s$. We now need some non-trivial
estimate towards the Ramanujan conjecture. For example, by a recent
result of Kim and Shahidi (\cite{Kim-Shahidi}), we know that
$$
- {1 \over 6} \leq \Re(s) \leq {1 \over 6}.
$$
This implies that
\begin{equation*}
\left| \varphi_v \begin{pmatrix} 1 \\ & \varpi_v \end{pmatrix}
\right| \leq {2q^{-{1 \over 2}+{1 \over 6}} \over 1 + q^{-1}} <
2q^{-{1 \over 3}}.
\end{equation*}
When $\G = \PGL_2$ and $\pi$ is cuspidal, the result follows from
the same result of Kim and Shahidi. Next, let $\G$ be an inner form
of the quasi-split group $\G'={\rm PGU}(2,1)$.
Suppose that $\G'$ splits over a quadratic
extension $E/F$. By Rogawski's theorem \cite{Rogawski}, there is
an automorphic cuspidal representation $\pi'= \otimes_v \pi_v'$ of
${\rm PGU(2, 1)}$ such that for all $v \notin S$, $\pi_v = \pi'_v$
as representations of $\G'(F_v)$. Consider the base change of
$\pi'$ from ${\rm PGU}(2,1)/F$ to $\PGL_3/E$, again established by
Rogawski. At this point, use the result of Oh on $\PGL_3$. Notice that this
is sufficient as the groups considered above at the only groups that
are of rank one at infinitely many places. \\

{\em The case where the s.s. rank is larger than one. } By
Proposition~\ref{one-dimensional}, $\pi_v$ is not one-dimensional
(for $v\notin S_F$), and the assertion follows from
Corollary~\ref{n}.
\end{proof}

\part{Geometry and height functions}

\section{Geometry}
\label{sect:geometry}


In this section we recall the constructions and basic geometric
properties of wonderful compactifications.
%
%

\subsection{Flag varieties}
\label{sect:flag}

An important class of varieties, homogeneous for the action of $\G$
is the class of generalized flag varieties
$$
Y_I:=\P_I\ba \G.
$$
The geometry of these, and their subvarieties (for example, Schubert
varieties) plays an important role in different branches of algebra,
e.g., representation theory and enumerative geometry.

\

We now recall some basic facts about these varieties. For $\la\in
\X^*(\T^{sc})$  we can define a line bundle $L_\la$ on
$\B^{sc}\ba\G^{sc}=\B\ba \G$ by
$$
\G^{sc}\times \ga /\sim, \,\,\, \text{ with } \,\,\, (g,a)\sim (gb,
\la^{-1}(b)a),
$$
$g\in \G^{sc}, b\in \B^{sc}$ and $a\in \ga$. The canonical
projection
$$
\pi\,:\, L_\la\ra \B\ba \G
$$
is given by
$$
\pi(g,a) = \B g.
$$
This gives an identification of
$$
\Pic(Y)=\X^*(\T^{sc}).
$$
Under this identification, the (closures of the) ample and the
effective cones of $Y$ correspond to the positive Weyl chamber, that
is, the set of nonnegative linear combinations of $\omega_i$. The
anticanonical class is given by
$$
-K_Y = 2\rho.
$$

\subsection{Wonderful compactifications}
\label{sect:wond}

\

First we work over an algebraically closed field
of characteristic zero.

\begin{prop}\cite{concini-p83},\cite{brion93}
\label{prop:geom} There exists a canonical compactification of a
connected adjoint group $\G$: a smooth projective variety $X$
such that
\begin{itemize}
\item $\G\subset X$ is a  Zariski open subvariety and the action of
$\G\times \G $ on $\G$, by
$$
(g_1,g_2)(g) =g_1gg_2^{-1},
$$
extends to
an action of $\G\times \G$ on $X$; \item The boundary $X\setminus
\G$ is a union of strict normal crossings divisors $D_i$ (for
$i=1,...,\sr$). For every $I\subset [1,...,\sr]$ the subvariety
$D_I=\cap_{i\in I}D_i$ is a  $\G\times \G$-orbit closure. All
$\G\times \G$-orbit closures are obtained this way;
\item $X$ contains a unique closed $\G\times \G$-orbit
$Y=\G/\B\times \G/\B$; \item The components $D_I$ are isomorphic to
fibrations over $\G/\P_I\times \G/\P_I$ with fibers canonical
compactifications of the adjoint form of the associated Levi groups.
\end{itemize}
\end{prop}

We review several constructions of wonderful compactifications over
algebraically closed fields (of characteristic 0), throughout $\G$
is semi-simple adjoint:

\

{\bf Via Hilbert schemes:} Let $\P\subset \G$ be a parabolic
subgroup and $Y=\G/\P$ the associated flag variety. Then $X$ is the
$\G\times \G$-orbit closure of the diagonal of the Hilbert scheme of
$Y\times Y$ (see \cite{brion97}).

\

{\bf Via representations:} Let $\lambda$ be a regular dominant
weight of $\G$ and $V_{\lambda}$ the irreducible representation of
$\G$ with highest weight $\lambda$. We have an action of $\G\times
\G$ on
$$
\End(V_{\lambda})=V_{\lambda}\otimes V_{\lambda}^*.
$$
Taking the closure of the orbit through (the image of) the identity
in $\mathbb P(\End(V_{\lambda}))$ we obtain the canonical
compactification $X$.

\

{\bf Via Lie algebras:} Let $\mg=\Lie(\G)$
and  $\sn=\dim(\G)=\dim(\mg)$. The variety $\mathbb L$ of
Lie subalgebras of the Lie algebra $\mg\oplus \mg$ can be regarded
as a subvariety of the Grassmannian ${\rm Gr}(n,2n)$. It contains
$\mg$, embedded diagonally. Moreover, ${\mathbb L}$ is a projective
$\G\times \G$-variety. Taking the closure of the $\G\times \G$-orbit
$X^\circ$ through $\mg$ we obtain $X\subset {\mathbb L}$. Since $\G$
is an adjoint group, the adjoint representation of $\G$ on $\G\times
\G$ is faithful and we may identify $X^\circ$ with the variety $\G$,
or more precisely $\G\times \G/{\rm Diag}(\G)$.

\

We proceed to describe the boundary $X\setminus \G$ in the latter
representation. Let $\P=\P_I$  (with $I\subset [1, \dots, r]$) be a
{\em standard} parabolic subgroup of $\G$. Choose  a Levi
decomposition $\P=\M\cdot \U$. Let $\mathfrak L_{P}$ be the set of
pairs
$$
(m+u,m+u'), \,\, \text{ with } \,\, m\in \Lie(\M), \,\,\, u,u'\in
\Lie(\U)  .
$$
Then  $\mathfrak L_{\P}$ is a subalgebra of $\mg\oplus \mg$. Next
let $\cC=\{ \P\} $ denote the conjugacy class of parabolic subgroups
of $\G$ containing $\P$. We note that the $\G\times \G$-orbit of
$\mathfrak L_{\P}$ in $\mathbb L$ does not depend on the particular
choice of $\M$ in the Levi decomposition of $\P$; in fact, it
depends only on  the class $\{\P\}$. We denote this orbit by
$D_{\cC}$. When $\P$ is maximal, the orbit $D_{\cC}$ is a smooth
irreducible divisor in $X$. Moreover,
$$
X\setminus \G=\bigcup_\cC D_{\cC} \ ,
$$
the union over classes of maximal parabolics. As the classes $\cC$
of maximal parabolics in $\G$ are in bijection with the simple roots
$\al$ of $\T$, we may write then
$$
D_{\cC}=D_{i},
$$
if $\cC$ is the maximal parabolic that corresponds to $\al_i$.

\label{sect:pic}

\

For $\sw\in \sW $ we denote by $X(\sw)$ the closure of
$\B\sw\B\subset \G$ in $X$. The $\B\times \B$ - stable boundary
components $D_i$ correspond to $X( s\sw_0 s_i)$. Every line bundle
$L$ on $X$ restricts to the unique closed $\G\times \G$-orbit
$Y=\G/\B\times \G/\B$; we get a restriction map $\Pic(X)\ra
\Pic(Y)$. Recall, that in Section~\ref{sect:flag} we have identified
the Picard group of $\G/\B$ with $\X^*(\T^{sc})$.

\begin{prop}
\label{prop:geometry} Let $X$ be the canonical compactification of
$\G$ as above.

\begin{itemize}
\item The image of $\Pic(X)\hookrightarrow \Pic(Y)$ consists of
classes
$$
L(\lambda)=(\lambda, -\sw_0\lambda) \subset
\Pic(Y)=\Pic(\G/\B)\times\Pic(\G/\B);
$$
\item The (closed) cone of effective divisors is given by
$$
\Lambda_{\rm eff}(X):=\oplus_{i=1}^\sr\R_{\ge 0} [D_i].
$$
More precisely, if $\la\in \X(\T^{sc})$ is a dominant weight then
the line bundle $L(\la)$ on $X$ has a unique (up to scalars) global
section $f_\la$ with divisor
$$
\dv(f_\la)=\sum_{i=1}^\sr \langle \la, \alpha_i^{\vee}\rangle D_i.
$$
Moreover, $f_{\la}$ is an eigenvector of $\B^{sc}\times \B^{sc}$
with weight $(-\sw_0\la,\la)$. \item The anticanonical class is
given by
$$
-K_X = L(2\rho+\sum_{i=1}^r\alpha_i).
$$
\end{itemize}
\end{prop}

\begin{proof}
See \cite{concini-p83} and \cite{brion93}.
\end{proof}

\label{sectLb}

\


Let $L$ be a line bundle on $X$. Then $L$
admits a $\G$-linearization - there exists a $\G^{sc}\times \G^{sc}$-action on
$\pi\,:\, L\rightarrow X$
such that for $(g_1,g_2)\in \G^{sc}\times \G^{sc}, l\in L$ one has
$$
\pi((g_1,g_2)\circ l)=(\bar{g}_1,\bar{g}_2)\circ \pi(l)  \ ,
$$
where $\bar{g}$ is the image of $g\in \G^{sc}$ in $\G$. We have a
representation $\rho=\rho_L$ of $\G^{sc}\times \G^{sc}$ on the space
of global sections $H^0(X,L)$.

\begin{thm}
\label{theo:lbb} Let $\la\in \mX^*(\T^{sc})$ and $L(\la)$ be the
associated line bundle on $X$. Let $\rho_{\la}$ be the
representation of $\G^{sc}\times \G^{sc}$ on $H^0(X,L(\la))$. Then:
\begin{itemize}
\item The representation $\rho_{\la}$ decomposes with multiplicity
one;
\item Let $\gamma$ be the dominant weight on $\T^{sc}$
(relative to $\B^{sc}$). Let $\xi_{\gamma}$ be the associated
irreducible representation of $\G^{sc}$ and $\eta_{\gamma}$ the
irreducible representation of $\G^{sc}\times \G^{sc}$ defined by
$$
\eta_{\gamma}=\xi_{\gamma}\times \xi^*_{\gamma} \ .
$$
Then each irreducible component of $\rho_{\la}$ is of the form
$\eta_{\gamma}$ for some dominant weight $\gamma$ of $\T^{sc}$.
Moreover, if $\eta_{\gamma}$ appears in $\rho_{\la}$, then $\gamma$
has the form
$$
\gamma=\la-\sum_{i=1}^r n_i \la_i \ ,
$$
with $n_i\in \Z_{\ge 0}$, for all $i=1,...,r$.
\item In particular,
the restriction of $\rho_{\la}$ to $\G^{sc}\times \{1\}$ is a sum,
with multiplicity, of $\xi_{\gamma}$'s with $\gamma$ as above.
\end{itemize}
\end{thm}


\label{sect:non-split}

\

We now return to the case of nonsplit $\G$ and
discuss the canonical $F$-structure on
the wonderful compactification $X$ of $\G$ as well as the corresponding
Galois action on the boundary divisors $D_{\al}$, for $\alpha \in
\Delta$.


\

A summary can be given as follows: in the Lie
algebra model, the Galois group $\Ga$ operates on $\Gr(n,2n)$ in
an obvious way and preserves $\mathbb L$ and $X^{\circ}$. Hence
$\Ga$ operates on $X$. Since $\Ga$ acts on the parabolic
subgroups of $\G$, it also permutes the boundary divisors
$$
\sigma(D_{\al})=D_{\sigma(\al)}.
$$
Here if $\cC=\{\P\}$, then  $\sigma(\cC)=\{ \sigma(\P)\}$.

\

The $F$-group $\G'$ is defined by a homomorphism
$$
\theta\,: \Gamma \ra \Out(\G^{sp}) \ ,
$$
to the group of outer automorphisms relative to the pair
$(\B^{sp},\T^{sp})$ and $\B^{sp},\T^{sp}$ chosen as above and both
defined over $F$. We have then an action of $\Gamma$ on $\Delta$.
Moreover,
$$
\sigma(D_{\al})=D_{\sigma({\al})} \,\, \text{ for } \al \in \Delta \
.
$$

\ Finally, since $\G$ is obtained $\G'$ by inner twisting, the two
actions of $\Gamma$ on the classes of maximal parabolics in
$\G(E)$ and $\G'(E)$ coincide.

\

In detail, let $\la=\sum_{i=1}^{\sr} \omega_i$ and $V_{\la}$ be the
corresponding representation with highest weight $\la$. We first
regard $V_{\la}$ as a module for $\G'$. As such, it is defined over
$F$, since $\la\in \mathfrak X^*(\T')_F$ (the action of $\Ga$ on
$\mathfrak X^*(\T')=\mathfrak X^*(\T)$ via $\theta$ simply permutes
the $\omega_i$). Let $L_{2\la}$ be the corresponding line bundle  on
$\B'\ba \G'$. Then $W:=H^0(\B'\ba \G', L_{2\la})$ has an
$F$-structure $W_F$. Denote by $\rho=\rho_{\la}$ the absolute
representation of $\G$ on $W$. By the Borel-Weil theorem, we may
identify $V_{2\la}$ with $W_F\otimes_F \bar{F}$, this gives an
action of $\Gamma$ on $V_{2\la}$:
$$
\sigma(w\otimes a) = w\otimes \sigma(a), \,\, \text{ for } w\in W_F,
a\in \bar{F}.
$$
We have the twisted action $\tilde{\sigma}=c(\sigma)\cdot \sigma$ of
$\Gamma$ on $\G(\bar{F})$, the corresponding twisted action on
$W\otimes W^*$ is given by
$$
\tilde{\sigma}\cdot w =\rho\otimes \rho^*(c(\sigma),c(\sigma))\cdot
\sigma.
$$
If we identify $\End(W)=W\otimes W^*$ we see that this action is
\begin{equation}
\label{eqn:8.1} \tilde{\sigma}(A)
=\rho(c(\sigma))\sigma(A)\rho(c(\sigma))^{-1}.
\end{equation}
Let $e_W\in \End(W)$ be the identity. Since $W$ is defined over $F$,
$\sigma(e_W)=e_W$ and similarly, $\tilde{\sigma}(e_W)e_W$. Thus
$e_W$ is rational for the twisted action.

\begin{rem}
The representation $\rho$ is actually a representation of $\G^{sc}$.
However, $\rho$ is irreducible and therefore maps $\rZ^{sc}$ to
scalars. It follows that \eqref{eqn:8.1} is well-defined.
\end{rem}

We use the twisted action to define an $F$-structure on $\End(W)$.
We claim then that the map
$$
\G\ni g\mapsto \rho(g)\cdot e_W=\rho(g)
$$
is $F$-rational. Consider
\begin{align*}
\rho(c(\sigma)\sigma(g))e_W
 = &  \rho(c(\sigma)\sigma(g) c(\sigma)^{-1})e_W  \\
 = &  \rho(c(\sigma))\rho(\sigma(g))e_W\rho(c(\sigma))^{-1} \\
 = &  \rho(c(\sigma))\sigma(\rho(g)e_W)\rho(c(\sigma))^{-1}.
\end{align*}
The latter equality holds since $e_W$ is $F$-rational. Thus
$$
\rho(c(\sigma)\sigma(g))e_W = \tilde{\sigma}(\rho(g)e_W).
$$
which proves the assertion.

For $g\in \G$, let $\bar{\rho}(g)$ be the image of $\rho(g)$ in
$\mathbb P(W)$. Then
$$
\bar{\rho}\,:\, \G\ra \mathbb P(\End(W))
$$
is also $F$-rational (equivalently, $\bar{\rho}$ is
$\Ga$-equivariant for the twisted action). Thus $\bar{\rho}(\G)$ and
its closure $X$ in $\mathbb P(\End(W))$ are $\Ga$-stable. This
gives a canonical $F$-structure for the variety $X$ and on
$$
Y:=\bigcap_{\al\in \Delta} D_{\al},
$$
(the $\G\times\G$-orbit through the image of $v_{2\la}\otimes
v_{2\la}^*$ in $\mathbb P(\End(W))$, where $v_{2\la}$ is ``the''
highest weight in $V_{2\la}$, see 5.1 of \cite{concini-p83}).
Moreover, we have an equivariant isomorphism of $F$-varieties
$$
Y\simeq\G/\B\times \G/\B^{-}.
$$

Indeed, the $\bar{F}$-irreducible boundary components $D_{\al}$ of
$X\setminus \G$ are permuted by $\Gamma$. Each component $D_{\al}$
is $\G\times \G$-stable and in particular stable for the action of
the diagonal ${\rm Diag}(\G)$:
$$
A\mapsto \rho(g)A\rho^{-1}(g).
$$
From \eqref{eqn:8.1} we have
$$
\tilde{\sigma}(D_{\al})= \sigma(D_{\al}).
$$
Thus $Y$ is invariant for the {\em standard} action.

We now prove that the $\G\times \G$-action on $X$ is $F$-rational
for the twisted $F$-structure. For this we need to see that for
$x\in X$,
$$
\tilde{\sigma}((g_1,g_2))x=\tilde{\sigma}(g_1,g_2)\tilde{\sigma}(x).
$$
In fact, this holds for all $y\in \mathbb P(\End(W))$. We have to
show that
$$
\tilde{\sigma}(\rho(g_1)A\rho(g_2)^{-1}) =
\rho(\tilde{\sigma}(g_1))\tilde{\sigma}(A)\rho(\tilde{\sigma}(g_2)^{-1}).
$$
The left side is
$$
\rho(c(\sigma))\sigma(\rho(g_1)A\rho(g_2)^{-1})\rho(c(\sigma))^{-1},
$$
and the right side
$$
\rho(c(\sigma))\rho(\sigma(g_1))\rho(c(\sigma))^{-1} \rho
c(\sigma)(A)\rho(c(\sigma))^{-1}
\rho(c(\sigma))\rho(\sigma(g_2))^{-1}\rho(c(\sigma))^{-1} =
$$
$$
\rho(c(\sigma))\rho(\sigma(g_1))\sigma(A)\rho(\sigma(g_2))^{-1}\rho(c(\sigma))^{-1}.
$$
The assertion follows, since $\sigma $ and $\rho(g)$, $g\in \G$,
commute.

\

The restriction of line bundles to the unique closed $\G\times
\G$-orbit $Y$ induces an injection
\begin{equation}\label{eqn:inj}
\Pic(X)\ra \Pic(Y).
\end{equation}
Note that since each divisor $D_{\al}$ is $\G\times \G$-stable, the
two actions of $\Ga$ on $X$ give rise to the same action on
$\Pic(X)$. Moreover, since $Y$ is defined over $F$, the injection
\eqref{eqn:inj} is $\Ga$-equivariant. Here
$$
\Pic(Y)\simeq \mathfrak X(\T')\oplus \mathfrak X(\T')
$$
as $\Ga$-modules. The image of $\Pic(X)$ is exactly the set of
pairs $(\la,-\la)$, where $\la\in \mathfrak X^*(\T')$, with boundary
divisors $D_{\al}$ corresponding to $(\al,-\al)$. The
$F$-irreducible boundary components are divisors
$$
D_J=\sum_{\al\in J}D_{\al},
$$
for any $\Ga$-stable subset $J\subset \Delta(\G',\T')$.

\section{Heights}
\label{sect:heights}

\subsection{Metrizations}
\label{sect:metr}

Here we recall the definitions of (adelically) metrized line
bundles and the associated heights.

\begin{defn}
\label{defn:metri} Let $X$ be a smooth projective algebraic variety
over a number field $F$. A smooth adelic metrization of a line
bundle $L$ on $X$ is a family of $v$-adic norms $\|\cdot \|_v$ on
$L\otimes_F{F_v}$ for all $v\in \Val(F)$ such that

\begin{itemize}
\item for $v\in S_{\infty}$ one has $\|\cdot \|_v$ is
$C^{\infty}$;
\item for $v\notin S_{\infty}$ the norm
of any local section of $L$ is locally constant in $v$-adic
topology;
\item there exist a finite set $S\subset \Val(F)$,
a flat projective scheme (an integral model) $\mathcal X$ over
$\Spec(\cO_S)$ with generic fiber $X$ together with a line bundle
$\cL$ on ${\mathcal X}$ such that for all $v\notin S$ the $v$-adic
metric is given by the integral model.
\end{itemize}
\end{defn}

\begin{exam}
\label{exam:metr} If $L$ is generated by global sections $(\mathsf
s_i)$ and $\mathsf s$ is a section such that $\mathsf s(x)\neq 0$
then
$$
\|\mathsf {\mathsf s}(x)\|_v:=\max_{i}(|\frac{\mathsf s_i}{\mathsf
s}(x)|_v)^{-1}.
$$
This defines a $v$-adic metric on $L$, which, of course, depends on
the choice of the basis $(\mathsf s_i)$. An adelic metric on $L$
is a collection of $v$-adic metrics (for all $v$) such that there
exists an $F$-rational basis $(\mathsf s_j)$ of $H^0(X,L)$ with the
property that for {\em all but finitely many} $v$ the $v$-adic
metric on $L$ is defined by means of this basis.
\end{exam}

An adelically metrized line bundle $\cL$ induces local and
global heights: for any local section $\mathsf s$ of $L$ and any $x$
with $\mathsf s(x)\neq 0$ define
$$
H_{\fs,\cL,v}(x)=\|\mathsf s(x)\|_v^{-1}.
$$
For $x\in X(F)$ the product formula ensures that the global height
$$
H_{\cL}(x)=\prod_{v\in \Val(F)}H_{\fs,\cL,v}(x)
$$
is independent of the choice of $\fs$. We write $\cL=(L,\|\cdot\|)$
when we want to emphasize that $L\in \Pic(X)$ is adelically
metrized.

\

Let $F$ be a number field and $V$ a finite-dimensional vector space
over $F$. Thus $V$ is the set of $F$-rational points $V(F)$ of a
linear variety $V$ defined over $F$. For $v\in \Val(F)$ we set
$V_v:=V\otimes _F F_v$.

Suppose first that $v$ is non-archimedean. Let $\La_v$ be an $\cO_v$
lattice in $V_v$. We define the norm
$$
\|\cdot \|_{\La_v} =\| \cdot \|_v
$$
on $V_v$ (associated to $\La_v$) as follows. Let $n=\dim V$ and let
$$
\cB_v:=\{ \xi_v^{(1)},\ldots,\xi_v^{(n)}\}
$$
be an $\cO_v$-module basis of $\La_v$. Let $v\in V_v$ have the form
$$
v=\sum_{j=1}^n a_j\xi_{v}^{(j)} \ .
$$
We set
$$
\|v\|_{\La_v}=\max_{1\le j\le n} (|a_j|_{v})  \ .
$$
We see at once that $\|\cdot \|_{\La_v}$ depends only on $\La_{v}$
and not on the particular choice of a basis of $\cB_{v}$.

Now suppose that $v$ is archimedean and let
$$
\cB_{v}:=\{ \xi_v^{(1)},\ldots,\xi_v^{(n)}\}
$$
be a basis of the $F_{v}$-vector space $V_{v}$. For
$$
v=\sum_{j=1}^n a_j\xi_{v}^{(j)}
$$
we set
$$
\|v\|_{\La_v}=(\sum_{1\le j\le n} |a_j|^2_{v})^{1/2} \
$$
(for $v$ complex, we set $|z|^2_{v}=z\cdot \bar{z}$).

\subsection{Heights on the canonical compactification}
\label{sect:hhe}

First we describe the situation for $\G^{sp}$. Recall that the
(classes of) irreducible boundary components $D_{\al}$, $\al\in
\Delta=\Delta(\G^{sp},\T^{sp})$, generate $\Pic(X)$. Each divisor
class $[D_{\al}]$ has a unique $\G^{sp}\times \G^{sp}$-stable
representative - namely this component. We put $L_{\al}=\mathcal
O_X(D_{\al})$. This gives us a canonical splitting of the projection
$$
\Div(X)\ra \Pic(X),
$$
which allows us to identify $\Pic(X)$ with the set of all integral
linear combinations
$$
L=\sum_{\al\in \Delta} s_{\al}L_{\al},
$$
and $\Pic(X)_{\C}$ with formal sums $\sum_{\al\in \Delta} s_{\al}
L_{\al}$.

\

We fix an integral model $\mathcal L_{\al}$ on the line bundle
$L_{\al}$, for $\al\in \Delta$. If $\mathcal L_{\al}'$ is another
model then the induced integral structures on $\cL_{\al,v}$, resp.
$\cL'_{\al,v}$ coincide for almost all $v$. An integral model
defines a height function
$$
H_{\al}=H_{\cL_{\al}}\,:\, X(F) \ra  \R_{>0}
$$
as in Section~\ref{sect:metr}. Given an $L=\sum_{\al\in \Delta}
s_{\al}L_{\al}$, with $s_{\al}\in \C$, we may define a height
$$
\begin{array}{rcc}
H_{L}   \,:\, X(F)& \ra &  \C \\
               x  & \mapsto & \prod_{\al} H_{\al}(x)^{s_{\al}}
\end{array}
$$

\

We make the above construction more precise and explicit, by
defining local heights
$$
H_{\al,v}\,:\, \G(F_v)\ra \R_{>0}, \,\,\,\text{ for all } \,\,\, \al.
$$
For each $L_{\al}$ we fix the (unique, up to scalars) $F$-rational
global section $\fs_{\al}\in H^0(X,L_{\al})$, which is
$\G^{sp}\times \G^{sp}$-invariant and non-vanishing on $\G^{sp}$.
Using the integral structure, put, for $g_v\in \G^{sp}(F_v)$ and
$\al\in \Delta$,
$$
H_{\al,v}(g_v):=\|\fs_{\al}(g_v)\|_v^{-1}\,\,\text{ and }\,\,
H_{\al}:=\prod_v H_{\al,v}.
$$

\

\begin{prop}
\label{prop:compare} There is a lattice $\La$ in some representation
of $\G^{sp}$ with the following property: Let $\K_f$ be the
stabilizer of $\La_f$ in $\G^{sp}(\A_f)$. Then for almost all finite
$v\in \Val(F)$ and for every
$$L=\sum_{\al\in \Delta} n_{\al}L_\al\in \Pic(X), \,\,\, \text{ with } n_{\al}\ge 0
$$
we have
$$
H_{L,v} (k_1 gk_2)=H_{L,v}(g), \,\,\,\text{ for all } \,\, g\in
\G(F_v), \,\,\,k_1,k_2\in \K_v.
$$
\end{prop}

\begin{proof}
If suffices to consider $L=L_{\la}$, for regular dominant weights
$\la$. We may assume that the action of $\G^{sp}\times \G^{sp}$ on
$L$ is defined over $F$ (passing to a multiple of $L$, see
\cite{mumford}, Section 3, Prop. 1.5). Fix a lattice $\Lambda
\subset H^0(X,L)_F$, it defines an $\cO_F$-integral structure on
$L$.
\end{proof}

\subsection{Local heights}

Let $L=L_{\la}$ be a very ample $F$-rational line bundle on $X$ and
fix a lattice $\La\subset H^0(X,L)$. Over the (fixed) splitting
field $E$, there is a distinguished $\fs \in H^0(X,L)$ which is
$\G^{sc}\times \G^{sc}$-invariant and non-vanishing on $\G(E)$.
Moreover, $\fs$ can be written as a product of sections
$\fs_{\al}\in H^0(X,L_{\al})$, $\al\in \Delta$, with support in the
$E$-rational divisor in $D_{\al}$. Thus, for $\sigma\in
\Ga=\Ga_{E/F}$ we have
$$
\sigma(\fs)=c(\sigma) \cdot \fs
$$
for some $c(\sigma)\in E^\times$. By Hilbert's theorem 90, we may
assume, after replacing $\fs$ by a suitable multiple, that $\fs$ is
fixed by $\Ga$, i.e.,  $\fs$ is $F$-rational.

We have an $\cO_E$-integral structure on $L$ over $E$, induced from
$\La\otimes_{\cO_F} \cO_E$. Over $E$, we have a decomposition
$$
H^0(X,L_{\la})=\oplus_\gamma (V_{\gamma}\otimes V^*_{\gamma}),
$$
as $\G^{sc}\times \G^{sc}$-modules. Here the sum is over
dominant $\gamma$ of the form $\gamma=\la-\sum_{\al\in \Delta}
m_{\al} \al$, with $m_{\al} \ge 0$. The following lemma is used to
calculate the local height integrals.
\begin{lem}\label{lem:height-formula}
One may choose $S_F$ large enough, so that for $v \notin S_F$, if
$g_v=k_vt_vk_v'$ with $k_v \in \K_v$ and $t_v \in \S(F_v)^+$, then
$$
H_v(g_v)=|\chi_{\la}(t_v)|_w.
$$
Here $\chi_\la$ is the rational character of $\T$ associated with
the dominant weight $\lambda$.
\end{lem}

\begin{proof}
After replacing $L_{\la}$ by a positive integral multiple
we can assume that $\la$ is trivial on the center of $\G^{sc}$. In
particular, $\chi_\la(t_v)\in F_v^\times$. Since $L_{\la}$ is
$\G\times \G$-linearizable,  we have an $F$-rational representation
$\varrho$ of $\G\times \G$ on $H^0(X,L_{\la})$. Then for $g_v\in
\G(F_v)$ we have
$$
H_v(g_v)=\|\varrho(g_v)\|_{\La_v}
$$
where $\La_v=\La\otimes_{\cO_F} \cO_v$.

Let $\La_E\subset H^0(X,L)_E=V_E$ be an ``admissible'' lattice,
i.e., $\La_E$ is homogeneous with respect to the decomposition of
$V=\oplus_\gamma V_{\gamma}$ into weight spaces.  More precisely,
the representation of $\G^{sp}\times \G^{sp}$ on $V$ has a $\Z$-form
on $\La_{\Z}$ with $\La_{\Z}$ a homogeneous lattice in $V_{\Q}$. We
set $\La_{\cO_E}:= \La_\Z \otimes_{\Z}\cO_E$ and $\La_{\cO_F}
:=\La_{\cO_E}\cap V_F$. Let $\mathcal G$ be the group scheme
structure on the double stabilizer of $\La_{\cO_F}$ and $\K_v =
\mathcal G(\cO_v)$, $\K_w=\mathcal G(\cO_w)$. We choose $S_F$ so
that for all $v\notin S_F$ and $w \mid v$
$$
\La_w=\La_v\otimes_{\cO_v}\cO_w.
$$
We choose $|\cdot |_w$ so that its restriction to $F_v$ coincides
with $|\cdot |_v$. Then the restriction of $\|\cdot \|_{\La_w}$ to
$V_v$ is $\|\cdot \|_{\La_v}$ and, for $g\in \G(F_v)$
$$
H_v(g_v)=\|\varrho(g_v)\|_{\La_w}.
$$
Write, according to Corollary \ref{coro:comp},
$$
g_v=a_w x_v a_w^{-1}
$$
with $x_v\in \G'(F_v)$, $a_w\in \mathcal G(\cO_w)=\K_w$. Then
$$
H_v(g_v)=\|\varrho(a_w)\varrho(x_v)\varrho(a_w)^{-1}\|_{\La_w}=\|\varrho(x_v)\|_{\La_w}
$$
We set next
$$
\K'_v:=\K_w\cap \G'(F_v).
$$
Recall that $\G(F_v)\subset \G(E_w), \G'(F_v)\subset \G(E_w)$. By Bruhat-Tits theory, we have
$$
\G'(F_v)=\K'_v\S'(F_v)^+\K_v'.
$$
Write accordingly, $x_v=k_1t_vk_2$, with $k_1,k_2\in \K_v',\,\,
t_v\in \S'(F_v)^+$. Then
$$
H_v(g_v)=\|\varrho(t_v)\|_{\La_w}.
$$
Next choose a basis $\mathcal B=\{ v_{\mu}\}$ of $\La_w$ consisting
of weight vectors for $\T_v^{sp}$. Thus
$$
H_v(g_v)=\max_{\mu} |\mu(t_v)|_w
$$
Each $\mu$ has the form
$$
\mu=\chi_{\la}\cdot \prod_{\al\in \Delta}\al^{-m_{\al}},
$$
with $m_\al\in \N$. Since $t_v\in \S'_v(F_v)^{+}$, we have
$|\al(t_v)|_w=|\al(t_v)|_v\ge 1$. Hence finally
\begin{equation}
\label{eqn:loc-height} H_v(g_v)=|\chi_{\la}(t_v)|_v.
\end{equation} \end{proof}
\

\subsection{Local integrals representing heights} \label{sect:def-h} Let $\D$ be a
central simple algebra of rank $m$ over the number field $F$. Also
let $\La$ be an arbitrary lattice in $\D$. We set, for each place
$v$, $\D_v = \D \otimes_F F_v$, and if $v$ is non-archimedean,
$\La_v = \La \otimes_{\cO_{F}} \cO_v$. In this subsection, we define
a family of norms $ \| \cdot \|_{\Lambda_{v}}$ on $\D_v$, one for
each place $v$ of $F$, subject to a certain compatibility condition.

\begin{itemize}
\item {\em non-archimedean $v$:} Choose a basis $\{ \xi_1^v,
\dots, \xi_k^v\}$ for $\D_v$ with $\xi_i^v\in \Lambda_{v} = \Lambda
\otimes_\cO \cO_v$ for all $i$. For $g \in \D_v$, write $g = \sum_{i}
c_i(g) \xi_i^v$ and set

\begin{equation*}
\|g\|_v=\| g \|_{\Lambda_{v}}:= \max_{i=1, \dots, k} \{ |c_i(g)|_v
\}.
\end{equation*}
It is easy to see that this norm is right and left
$\Lambda_{v}^\times$-invariant and therefore independent of the
choice of the basis.

\item {\em archimedean $v$:} Fix a Banach space norm $\|\cdot
\|_v=\|\cdot\|_{\D_v}$ on the finite dimensional real (or complex)
vector space $\D_v = \D \otimes_F F_v$.

\end{itemize}

 Clearly, for $c\in F_v$ and $g \in \D_v$,  we have

\begin{equation*}
\|  cg \|_v = |c|_v \cdot \|g\|_v.
\end{equation*}
Consequently, for $c \in F$ and $g \in D$, we have

\begin{equation}\label{artin}
\prod_v \|cg \|_v = \prod_v \|g\|_v,
\end{equation}
by the product formula. This is the compatibility condition
mentioned above. Define a function $\Psi_v$  as follows:
$$
\Psi_v =\left\{  \begin{array}{lcl} \text{characteristic function of
} \Lambda_{v}, &
{\rm for} & v \notin S_{\infty},\\
\exp( - \pi \|\cdot \|_{v}^2) & {\rm for} & v\in S_{\infty}.
\end{array}
\right.
$$

\begin{lem}
\label{height-schwartz} For all $v\in \Val(F)$, all $g_v \in \D_v$
and all $s$ with $\Re(s)>0$ one has
\begin{equation*}
\int_{F_v^\times} \Psi_v(a g_v) |a|^s \, d^\times a =
\zeta_{F,v}(s)\cdot \|g_v\|_v^{-s}.
\end{equation*}
Here $\zeta_{F,v}(s)$ is the $v$-local Euler factor of the zeta
function of $F$.
\end{lem}

\noindent For $\Phi_v \in C_c^\infty(\D_v)$ and $\Re(s)>0$ we set
\begin{equation*}
H_v(s, g; \Phi_v) := \zeta_{F, v}(s)^{-1} \int_{F_v^\times} \Phi_v(a
g) \vert a \vert^{s} \, d^\times a,
\end{equation*}
where $g \in \D_v$. For $g =(g_v)_v \in \D(\A)$ define the {\em global
height function}:

\begin{equation*}
H(g) =\prod_{v\in \Val(F)} H_v(g)= \prod_{v\in \Val(F)} \|g_v \|_v.
\end{equation*}
Similarly, if $\Phi = \otimes_v \Phi_v$ is a global Schwartz-Bruhat
function on $\D(\A)$, we define
\begin{equation*}
H(s, g; \Phi) = \prod_{v\in \Val(F)} H_v(s, g_v; \Phi_v).
\end{equation*}
We extend the functional $H(s, g;\cdot )$ to $C_c^\infty(\D(\A))$ by
linearity. By the product formula, both $H(g)$ and $H(s, g; \Phi)$,
(for $\Phi \in C_c^\infty(\D(\A))$) are well-defined on the
projective group of $\D$.

\subsection{Complexified height function}
\label{sect:complexified-height}

\begin{notationnum}
\label{nota:t}
Let $\mathcal{T}$ be the set of all Galois invariant
${\bf s}=(s_\alpha)_{\alpha \in \Delta(\G', \T')}$ (c.f.
\ref{sect:expression-infinite} for the definition of Galois
invariance). The element $\underline{0} \in \mathcal{T}$ is
defined by setting all coordinates equal to zero. For $\epsilon\in \R$,
the set $\mathcal{T}_\epsilon$ is the
set of ${\bf s} = (s_\alpha)_\al \in \mathcal{T}$ such that
$\Re(s_\alpha)
> \kappa_{\alpha} + 1+\epsilon$, for all $\alpha$.
Starting with an element ${\bf s} = (s_\alpha)_{\alpha \in \Delta(\G',
\T')}$ and $v \notin S$, we obtain a tuple ${\bf s}^v =
(s^v_\vartheta)$ indexed by $\Delta(\G'(F_v), \sS'_v(F_v))$ by
setting $s^v_{r_v (\iota^*(\alpha))} = s_\alpha$; this is
well-defined. For ${\bf s}, {\bf t} \in \mathcal{T}$ and $v \notin
S$, we set
$$
<{\bf s}, {\bf t}>_v = \sum_{\vartheta \in
\Delta(\G'(F_v), \S'_v(F_v))} s^v_\vartheta t^v_\vartheta.
$$
When
there is no danger of confusion, we write ${\bf s}$ for ${\bf s}^v$. For
each subset $R$ of $\C$, we set $\mathcal{T}(R)$ to be the
collection of ${\bf s}= (s_\alpha)_\alpha$ with $s_\alpha \in R$ for
all $\alpha$.
\end{notationnum}

We go back to Lemma \ref{lem:height-formula}. For $v\notin S_F$,
we have expressed $H_v(g) $ in terms of $\chi_{\la}\in \mathfrak
X^*(\T)$. Write
$$
\chi_{\la} = \prod_{\al\in \Delta(\G',\T')} \al^{n_{\al}}
$$
and $\la=\sum_{\al\in \Delta(\G,\T)} n_{\al}\al$, $n_{\al}\in
\mathbb N$. We know that $(n_\alpha)_{\alpha \in \Delta(\G', \T')} \in
\mathcal{T}$. Now, let ${\bf s}: =(s_\alpha)_{\alpha \in \Delta(\G',
\T')} \in \mathcal{T}$. Suppose $v \notin S_F$, and $g_v \in
\G'(F_v)$ is written as $k_v t_v k_v'$ with $t_v \in \S'_v(F_v)^+$
and $k_v, k_v' \in \K'_v$. We define
\begin{equation}
H_v({\bf s}, g_v) = \prod_{\alpha \in \Delta(\G', \T')}
|\alpha(t_v)|_v^{s_\alpha}.
\end{equation}
Observe that
\begin{equation}
H_v({\bf s}, g_v) = \prod_{\vartheta \in \Delta(\G'(F_v), \S_v'(F_v))}
|\vartheta(t_v)|_v^{\ell_v(\vartheta) s_\vartheta}.
\end{equation}
If we use Lemma~\ref{coro:comp} to identify the groups $\G(F_v)$
and $\G'(F_v)$, we get a complexified local height function on
$\G(F_v)$. Let $S = S_F$ and suppose $g = (g_v)_{v \notin S} \in
\G (\A_S)$. We define
\begin{equation}
H_S({\bf s}, g) : = \prod_{v \notin S} H_v({\bf s}, g_v).
\end{equation}

\

\subsection{Local height integrals I}
\label{sect:rk}
Consider the
integral of the complexified local height function
\begin{equation}
\mathcal{J}_v(s) := \int_{\G_v} H_v(\bfs, g)^{-1} \, dg
\end{equation}
and its versions.
\begin{thm}
\label{ramified} \begin{enumerate} \item For all $v\notin
S_{\infty}$ the integral $ \mathcal{J}_v(s)$ is a holomorphic
function of $\bfs$ for $\bfs \in \cT_{-1}$.
\item Let $v$ be an archimedean valuation and $\partial$ any
element of the universal enveloping algebra. Then
\begin{equation*}
\mathcal{J}_{v,\partial}(\bfs):= \int_{\G(F_v)} \partial( H(\bfs,
g)^{-1} )\, d g
\end{equation*}
is holomorphic for $\bfs \in \cT_{-1}$.
\end{enumerate}
\end{thm}
\begin{proof} We will only prove the first part; the second part is similar.
Locally, every two local integral structures give rise to
essentially equivalent height functions; so, we replace the local
integral structure so that the resulting height function is
invariant under $\sK_v$, a good maximal compact subgroup. Let
$\underline{\sigma}$ be the vector consisting of the real parts of
the components of $\bfs$. The local height integral is majorized by
\begin{align*}
& \sum_{t \in \sS(F_v)^+} \sum_{\omega \in \Omega_v}
H(\underline{\sigma}, t \omega )^{-1} \text{vol }(\sK_v t
\omega \sK_v) \\
& \ll  \sum_{t \in \sS(F_v)^+} H(\underline{\sigma}, t)^{-1} \de_B(t) \\
& = \prod_{\vartheta \in \Delta(\G(F_v),
\sS_v(F_v))}\sum_{l=0}^\infty
\delta_B(\check\vartheta (\varpi_v^l))H(\underline{\sigma}, \check\vartheta(\varpi_v^l))^{-1}\\
& =\prod_{\vartheta \in \Delta(\G(F_v),
\sS_v(F_v))}\sum_{l=0}^\infty q_v^{-(\sigma_\vartheta -
\kappa_\vartheta)l \ell(\vartheta)}.
\end{align*}
The result is now immediate. \end{proof}
\begin{coro}
In the non-archimedean situation, for each $\epsilon > 0$ there is a
constant $C_v(\eps)$, such that $|\mathcal{J}_v(\bfs)| \leq
C_v(\eps)$ for all $\bfs \in \cT_{-1+\eps}$. In the archimedean
situation, for all $\eps>0$ and all $\partial$ as above, there is a
constant $C_v(\partial, \eps)$ such that
$|\mathcal{J}_{v,\partial}(\bfs)| \leq C_v(\partial, \eps)$ for all
$\bfs \in \cT_{-1 + \eps}$.
\end{coro}

\subsection{The integral of the local height function II}
Let $\G$ be a connected algebraic group over $F_v$, $v$
outside a finite set of places, and $\mathfrak{g}$ its Lie
algebra of invariant vector fields. Let $X$ be a smooth equivariant
compactification of $\G$. Denote by $D = X \backslash \G$ the
boundary. We assume that $D$ is a divisor with strict normal
crossings. Let $\mathcal{T}_X$ be the tangent bundle of $X$. We have
a restriction map
\begin{equation}
H^0 (X, \mathcal{T}_X) \to \mathcal{T}_{X, 1} = \mathfrak{g}
\end{equation}
obtained by evaluating a vector field at the neutral element $1 \in
\G$. Conversely, given $\partial \in \mathfrak{g}$, there is a
unique vector field $\partial^X$ such that for any open subset $U$
of $X$ and for any $f \in \mathcal{O}_X(U)$, $\partial^X(f)(x) =
\partial_g f(g.x) |_{g=1}$. The map $\partial \mapsto \partial^X$
is a section of the restriction map. Let $\partial_1, \dots,
\partial_n$ be a basis for $\mathfrak{g}$. Then $\delta: =
\partial_1^X \wedge \dots \wedge \partial_n^X$ is a global section
of the line bundle $\det \mathcal{T}_X = K_X^{-1}$. Moreover,
$\delta$ does not vanish on $\G$. Because of these considerations,
if we know that $K_X^{-1}$ is ample, Peyre's Tamagawa measure
restricts to Weil's Tamagawa measure on $\G(F_v)$.

If $\mu_v$ is the local Weil-Tamagawa measure, for $v$ outside of a
finite set of places, we have
\begin{equation}
\mu_v(\G(\mathcal{O}_v)) = \frac{\# \G(k_v)}{q_v^{\dim \G}}.
\end{equation}
Therefore, we normalize the local measure by the appropriate factor
to guarantee that $\mu(\G(\mathcal{O}_v))=1$ for almost all $v$.

Let $A$ be a subset of $\Delta(\G, \T)$, we set
$$D_A = \bigcap_{\alpha \in A} D_\alpha$$
and
$$
D_A^0 = D_A \backslash \left(\bigcup_{A \subsetneqq A'}
D_{A'}\right),
$$
for $D_\alpha$ as defined in \ref{sect:non-split}. The
following theorem is the analog of Theorem 9.1. of
\cite{chambert-t02} in this situation, with the same proof
(see also Theorem 3.1. of \cite{Denef}).

\begin{thm}\label{motivic-integral}
We have
\begin{equation}\label{eq:motivic-integral}
\int_{\G(F_v)} H_v(\fs, g_v)^{-1}\, dg_v = \frac{1}{\# \G(\k_v)} \sum_{A}
\# D_A^0 (\k_v) \prod_{\alpha \in A} \frac{q_v - 1}{q_v^{s_\alpha -
\kappa_\alpha +1} - 1}.
\end{equation}
\end{thm}
\begin{proof}
We split the integral along residue classes modulo $\mathfrak{p}_v$.
Let $\tilde{x} \in X(\k_v)$, $\k_v$ the residue field of $F_v$, and $A
= \{ \alpha \, | \, \tilde{x} \in D_\alpha\}$, so that $\tilde{x}
\in D_A^0$.

We can introduce local (\'{e}tale) coordinates $x_\alpha$ ($\alpha
\in A$), and $y_\beta$ ($\beta \in B$) with $\# A + \# B = \dim X$
around $\tilde{x}$ such that, locally, the divisor $D_\alpha$ is
defined by the vanishing of $x_\alpha$. Then the local Tamagawa
measure identifies with the measure $\prod dx_\alpha \times \prod
dy_\beta$ on $\mathfrak{p}_v^A \times \mathfrak{p}_v^B$. If $dx$
denotes the fixed measure on $\G(F_v)$, one has the equality of
measures on $\G(F_v) \cap \, red^{-1}(\tilde{x})$:
\begin{equation}
dx = \frac{q_v^{\dim X}}{\# \G(\k_v)} H_v(\rho,x)d\mu_v =
\frac{q_v^{\dim X}}{\# \G(\k_v)} \prod_{\alpha \in A}
q_v^{(\kappa_\alpha + 1)v(x_\alpha)} \prod dx_\alpha \prod dy_\beta.
\end{equation}
Consequently,
\begin{align*}
\int_{red^{-1}(\tilde{x})} & H_v({\bf s},x)^{-1} d x \\
& = \frac{q_v^{\dim X}}{\# G(\k_v)}\int_{\mathfrak{p}_v^A \times
\mathfrak{p}_v^B} q_v^{-\sum_{\alpha \in A}(s_\alpha -\kappa_\alpha
- 1)v(x_\alpha) }\prod dx_\alpha
\prod dy_\beta \\
& =\frac{q_v^{\dim X}}{q_v^{\# B }\# G(\k_v)}\prod_{\alpha \in A}
\int_{\mathfrak{p}_v} q_v^{-(s_\alpha -\kappa_\alpha -
1)v(x_\alpha) }dx_\alpha \\
& =\frac{1}{\# G(\k_v)}\prod_{\alpha \in A}\frac{q_v -
1}{q_v^{s_\alpha - \kappa_\alpha} - 1 }.
\end{align*}

We have used the identity
\begin{align*}
\int_{\mathfrak{p}_v} q_v^{-s v(x)}\, dx & = \sum_{n=1}^\infty
q_v^{-sn} \vol \left( \mathfrak{p}_v^n \backslash
\mathfrak{p}_v^{n+1} \right) \\
& = \sum_{n=1}^\infty q_v^{-sn}q_v^{-n} \left( 1 - \frac{1}{q_v}
\right) \\
& = \frac{1}{q_v} \frac{q_v - 1}{q_v^{1+s} - 1 }.
\end{align*}
\end{proof}

\subsection{An application to volumes}
\label{sect:vol}

We start with the following lemma in the non-archimedean situation:
\begin{lem}
Let $\sK_v$ be a be a good maximal compact subgroup of $(\G(F_v),
\sS(F_v))$ so that the Cartan decomposition $\G = \sK_v \sS(F_v)^+
\Omega \sK_v$ holds. Normalize measures so that ${\text vol} (\sK_v)
=1$. Then for all $ad \in \sS(F_v)^+ \Omega$, we have
\begin{equation}
{\text vol}\, ( \sK_v ad \sK_v) \ll \delta_\B(a).
\end{equation}
\end{lem}

As $\Omega$ is a finite set, this is an immediate consequence of
Lemma 4.1.1 of \cite{silberger}. Below we need more detailed
information on the behavior of the the above volume in the
quasi-split situation. The following lemma suffices for our
purposes; the analogous statement for simply-connected is classical.
For split groups, the lemma is contained in \cite{Gross}.

\begin{lem}
There exists a constant $c$, independent of $v$, such that for all
$t \in \sS (F_v)^+$, one has
\begin{equation*}
\vol(\sK_v t \sK_v) \leq \delta_B(t)(1 + \frac{c}{q_v}).
\end{equation*}
\end{lem}
\begin{proof}We use Theorem \ref{motivic-integral}.
The left hand side of \eqref{eq:motivic-integral} is obviously equal
to
$$
\sum_{t \in \sS (F_v)^+} H(\fs, t)^{-1} \vol(\sK_v t \sK_v).
$$
The comparison of this expression with the right hand side of
\eqref{eq:motivic-integral} will give an explicit formula for the
volume from which our result will easily follow. Clearly the right
hand side of \eqref{eq:motivic-integral} is equal to
\begin{equation}
\sum_{t \in \S(F_v)^+} \frac{\# D_{A(t)}^0(\k_v)(q-1)^{\# A(t)}}{\#
\G(\k_v)}\delta_\B(t)H(\fs, t)^{-1}.
\end{equation}
Here $A(t) = \{ \delta \in \Delta(\G, \S); \delta(t) = 1\}$; we will
suppress dependence on $t$ and simply write $A$. Comparison gives
\begin{equation}
 \vol(\sK_v t \sK_v) =\frac{\# D_{A}^0(\k_v)(q-1)^{\# A(t)}}{\#
\G(k_v)}\delta_\B(t)
\end{equation}

The algebraic set $D_A$ is the fiber variety $D_A \to \G/ \P_A
\times \G/\P_A$, with fibers isomorphic to $\overline{\M}_A$, the
wonderful compactification of the adjoint group of $\M_A$, i.e.
$\M_A$ modulo its center. Here $\M_A$ is the Levi factor of the
parabolic subgroup $\P_A$. Then $D_A^0$ has fibers $\M_A$ modulo its
center. Consequently,
$$
\# D_A^0(\k_v) = (q-1)^{-\# A}\left[ \G(\k_v): \P_A(\k_v)\right]^2\cdot \#
\M_A(k_v).
$$
If $\U_A$ is the unipotent radical of $\P_A$, we have
$$
\frac{\# D_A^0 (\k_v)(q-1)^{\# A}}{\# \G(\k_v)} = \left[\G(\k_v):
\P_A(\k_v) \right] (\# \U_A(\k_v))^{-1}.
$$

Let $\sW_A$ be the Weyl group of the Levi factor $\M_A$. By the Bruhat
decomposition we have
\begin{equation}
\G(\k_v) = \bigcup_{w \in \sW_A \backslash \sW / \sW_A} \P_A(\k_v) w \left(
\U_A(\k_v) \cap w^{-1} \overline{\U}_A(\k_v)w\right),
\end{equation}
where $\overline{\U}_A$ the unipotent radical opposite to $\U_A$.
It follows that
\begin{equation}
\frac{\left[\G(\k_v): \P_A(\k_v) \right]}{\# \U_A(\k_v)} = \sum_{w \in
\sW_A \backslash \sW / \sW_A} \frac{\# \left( \U_A(\k_v) \cap w^{-1}
\overline{\U}_A(\k_v)w\right)}{\# \U_A(\k_v)}.
\end{equation}
As there is only one double coset of maximal dimension, we have
\begin{equation}
\frac{\left[\G(\k_v): \P_A(\k_v) \right]}{\# \U_A(\k_v)} \leq 1 +
\frac{\# \left(\sW_A \backslash \sW / \sW_A\right)}{q_v} .
\end{equation}
The lemma is now clear for split groups. In the quasi-split case, we
need only consider $\Ga$-stable subsets $A$. For
the complex parameters $s_\alpha$ we have the extra assumption that
$s_{\alpha} = s_{\sigma \alpha}$, for $\sigma \in \Ga$. The
proof of the lemma in the quasi-split case is similar,
and we omit it.
\end{proof}

\section{Regularization}
\label{sect:regularintegral}

In this section, for $v \notin S_F$,
with a (slight) abuse of notation we identify $\G(F_v)$ and
$\G'(F_v)$; this is permissible,
in light of Corollary \ref{coro:comp}.

\subsection{Integrals from one-dimensional representations}
\label{sect:one-dim-reg}
Let $\chi$ be a one-dimensional automorphic representation of $\G$, and $S = S_F$.
We proceed to study analytic
properties of the integral
\begin{equation}
\mathcal{J}_S({\bf s}, \chi) = \int_{\G(\A_S)} H_S({\bf s}, g)^{-1}
\chi(g) \, dg
\end{equation}
for Galois invariant ${\bf s}$.

\begin{thm}\label{thm:one-dim}
The product
\begin{equation}
\prod_{\mO} L(s_\mO, \xi_\mO(\chi))^{-1} \int_{\G(\A_S)} H_S({\bf s},
g)^{-1} \chi(g) \, dg
\end{equation}
is holomorphic on $\mathcal{T}_{-\epsilon}$ for some $\epsilon >
0$. The product is over all Galois orbits in $\Delta(\G', \T')$.
\end{thm}
\begin{proof}
We have
\begin{equation}
\mathcal{J}_S({\bf s}, \chi) = \prod_{v \notin S} \mathcal{J}_v({\bf s},
\chi_v)
\end{equation}
where
\begin{equation}
\mathcal{J}_v({\bf s}, \chi_v) = \int_{\G(F_v)} H_v({\bf s}, g_v)^{-1}
\chi_v(g_v) \, dg_v.
\end{equation}
Since $\G$ is of adjoint type, the collection of elements $\{
\check\vartheta(\varpi_v)\}_{\vartheta \in \Delta(\G(F_v),
\sS_v(F_v))}$ forms a basis for the semigroup
$\sS_v(F_v)^+$. For any vector ${\bf a} = (a_\alpha)_{\alpha} \in
\mathcal{T}(\N)$, we set
\begin{equation}
t_v({\bf a}) = \prod_{\vartheta \in \Delta(\G(F_v), \S_v(F_v))}
\check\vartheta(\varpi_v)^{a_\vartheta}.
\end{equation}
Write
\begin{equation*}
1+a_v = \sum_{{\bf a} \in \mathcal{T}(\N)} q_v ^{-\langle {\bf a},
{\bf s}\rangle_v} \chi_v(t_v({\bf a}))\delta(t_v({\bf a})),
\end{equation*}
and
\begin{equation*}
b_v = \sum_{{\bf a}\in \mathcal{T}(\N)} q_v ^{-\langle {\bf a},
{\bf s}\rangle_v} \chi_v(t_v({\bf a}))(\text{vol }(\sK_vt_v({\bf
a})\sK_v)- \delta(t_v({\bf a}))),
\end{equation*}
so that, by Cartan Decomposition, we have
\begin{equation}
\mathcal{J}_v({\bf s}, \chi_v) = 1 + a_v + b_v.
\end{equation}
Observe that
\begin{equation}\begin{split}\label{delta}
1 + a_v & = \sum_{\bf a} q_v^{-\langle {\bf a}, {\bf s}\rangle_v}
\chi_v(t_v({\bf a}))\delta(t_v({\bf a})) \\
& = \prod_\vartheta \biggl( \sum_{a_\vartheta =0}^\infty
\chi_v(\check\vartheta(\varpi_v))^{a_\vartheta}q_v ^{-(s_\vartheta - \kappa_{\vartheta}) a_\vartheta\ell(\vartheta)} \biggr) \\
& = \prod_\vartheta {1 \over 1-
\chi_v(\check\vartheta(\varpi_v))q_v^{-(s_\vartheta -
\kappa_{\vartheta})\ell(\vartheta)}}.
\end{split}\end{equation}
With Proposition \ref{Hecke-product} in mind, we proceed as follows.
Let $\mathsf{\sigma}= (\Re(s_\al))_\al$. Observe that in the
definition $b_v$ we may assume ${\bf a} \ne
\underline{0}$. Since for each $v \notin S$,
$$
\Biggl\{ {\bf a}\,\, | \,\, {\bf a} \ne \underline{0} \Biggr\} =
\bigcup_{\vartheta \in \Delta(\G'(F_v), \sS'_v(F_v))} \Biggl\{
{\bf a}; a^v_\vartheta \ne 0 \Biggr\},
$$
we have
\begin{align*}
\sum_{v \notin S} \Bigl| b_v \Bigr| & \leq \sum_{v \notin S}
\sum_\vartheta \sum_{a_\vartheta \ne 0} q_v ^{-\langle {\bf a},
{\bf \sigma}\rangle_v} \Biggl|
(\text{vol }(\sK_vt_v({\bf a})\sK_v)-  \delta(t_v({\bf a})))\Biggr|  \\
& \ll\sum_{ v \notin S} q_v^{-1} \sum_\vartheta \sum_{a_\vartheta
\ne 0} q_v ^{-\langle {\bf a}, {\bf \sigma}\rangle_v}
\delta(t_v({\bf a})) \quad{\text{ (by Section~\ref{sect:vol})}}\\
& =\sum_{v \notin S} q_v^{-1} \sum_\vartheta \biggl(
\sum_{a_\vartheta =1}^\infty q_v^{-(\sigma_\vartheta -
\kappa_{\vartheta})a_\vartheta\ell(\vartheta)} \biggr) \prod_{
\beta \ne \vartheta}
\biggl( \sum_{a_\beta=0}^\infty q_v^{-(\sigma_\beta - \kappa_{\beta})a_\beta\ell(\beta)} \biggr) \\
& =\sum_{v \notin S} q_v^{-1} \sum_\vartheta {
q_v^{-(\sigma_\vartheta - \kappa_{\vartheta})\ell(\vartheta)}
\over \prod_\beta
(1 -q_v^{-(\sigma_\beta - \kappa_{\beta})\ell(\beta)} ) } \\
& \ll \sum_\vartheta \sum_{v \notin S} q_v^{ - {3 \over 2}} < \infty.
\end{align*}
We need to show the existence of a $C>0$ such that $\vert 1+a_v
\vert \geq C >0$ for all $v$. For this
\begin{equation*}
\vert 1 + a_v \vert \geq  \prod_\vartheta {1 \over 1 + q_v^{-
\sigma_\vartheta + \kappa_{\vartheta}}} \geq \prod_\vartheta { 1
\over 2} \geq  {1 \over 2^{r}},
\end{equation*}
with $r = |\Delta(\G^{sp}, \T^{sp})|$. In fact, $q_v^{-\sigma_k +
\kappa_{\vartheta_k}} \leq q_v^{-{1 \over 2}} < 1$.

\

Note that for ${\bf s} \in \cT_{-\epsilon}$ the estimates are
uniform, i.e., the quotient
$$
\frac{\prod_{v \notin S} \mathcal{J}_v(\chi)}{ \prod_{v \notin S}
(1 + a_v)}
$$
is holomorphic in $\cT_{-\epsilon}$.
This finishes the proof of the theorem.
\end{proof}

\

\subsection{Integrals from infinite-dimensional
representations}
Let $\pi = \otimes_v \pi_v$ be an infinite-dimensional automorphic
representation of $\G$ which is not the automorphic representation
associated to an Eisenstein series when $\G = \PGL_2$
(this exceptional case was treated in \cite{PGL2}).
For $v \notin S_F$, let $\varphi_{v, \pi}$ be the normalized spherical
function associated to $\pi_v$ and set $S =
S_F$.
\begin{thm}
\label{non-triv} Let $c>0$ be as in Theorem \ref{bound}.
The infinite product
\begin{equation}\label{non-triv-eq}
\mathcal{J}_S({\bf s}, \pi) : = \prod_{v \notin S} \int_{\G(F_v)}
\varphi_{v, \pi}(g_v) H_v({\bf s}, g_v)^{-1}\, dg_v
\end{equation}
 is holomorphic for ${\bf s}\in \cT_{-c}$.
\end{thm}

\

\begin{proof} We use the notation of Section \ref{sect:one-dim-reg}. Set
\begin{equation}
\mathcal{J}_v({\bf s}, \pi) =\int_{\G(F_v)} \varphi_{v, \pi}(g_v)
H_v({\bf s}, g_v)^{-1}\, dg_v.
\end{equation}
By Cartan decomposition,
\begin{equation}
\mathcal{J}_v({\bf s}, \pi) = \sum_{{\bf a}}
\varphi_{v, \pi}(t_v({\bf a}))H_v({\bf s}, t_v({\bf a}))^{-1}\text{vol}(\K_v t_v(\bf a) \K_v).
\end{equation}

Let

\begin{equation*}
a_v = \sum_\vartheta q_v^{-s_\vartheta} \text{vol }( \sK_v
\check\vartheta(\varpi_v) \sK_v) \varphi_{v, \pi}
(\check\vartheta(\varpi_v))
\end{equation*}
and

\begin{equation*}
b_v = \mathcal{J}_v({\bf s}, \pi) - 1 - a_v.
\end{equation*}
Fix an $\epsilon < {1 \over 2}$. We claim that there is a set $S'$
of places of $F$ and a positive constant $C$ such that $\vert 1
+a_v \vert \geq C$ for all $v \notin S'$. Since $\vert 1 +a_v
\vert \geq 1 - \vert a_v \vert$, we only need to show that $\vert
a_v \vert$ is asymptotically bounded away from $1$. To see this we
use the fact that $\| \varphi_{v, \pi} \|_{\rL^\infty} \leq 1$.
Hence for $v$ outside a growing finite set $S'$

\begin{equation*}
\vert a_v \vert \leq \sum_\vartheta q_v^{-\sigma_\vartheta +
\kappa_{\vartheta}} \to 0,
\end{equation*}
for ${\bf s} \in {\cT}_{-\epsilon}$. This implies that one can choose
$S'$ such that if $v \notin S'$ and ${\bf s} \in {\cT}_{-\epsilon}$
then $\vert 1 + a_v \vert \geq C$, for a constant $C$ that depends
only on $\epsilon$.

Next, we prove that $\sum_v b_v$ is absolutely and uniformly
convergent on ${\cT}_{-\epsilon}$. Denote by $b(m)_v$ the set of
terms such that ${\bf a}^v$ has exactly $m$ nonzero coordinates.
We have

\begin{align*}
\sum_v \vert b(1)_v \vert & \leq \sum_\vartheta \sum_{v \notin S}
\sum_{a_\vartheta \geq 2}
q_v ^{-\sigma_\vartheta a_\vartheta + \kappa_{\vartheta} a_\vartheta},  \\
\intertext{Note that $a_\vartheta \ne 1$, because otherwise the
term would appear in $a_v$. This implies that } \sum_v \vert
b(1)_v \vert & \leq \sum_\vartheta \sum_{v \notin S} { q_v
^{-2\sigma_\vartheta + 2\kappa_{\vartheta}}
\over 1 - q_v^{-\sigma_\vartheta + \kappa_{\vartheta}}} \\
& \ll \sum_\vartheta \sum_{v \notin S}
q_v ^{-2\sigma_\vartheta + 2\kappa_{\vartheta}} \\
& \ll \sum_{v \notin S} q_v^{-2 + 2\epsilon} \\
& < \infty,
\end{align*}
for $\epsilon < {1 \over 2}$.

\

To continue, we fix an ordering of the
elements of $\Delta(\G(F_v), \sS_v(F_v))$, say $\{\vartheta_1,
\dots, \vartheta_l\}$. Next we verify the claim for $\sum_{m\ge
2}\sum_v \vert b(m)_v \vert$, where
\begin{equation*}
b(m)_v = \sum_{1 \leq i_1 < i_2 < \dots i_m \leq l} \sum_{a_{i_1}
>0} \sum_{a_{i_2}>0} \dots \sum_{a_{i_m}>0} q_v^{-\sum_{j=1}^m
a_{i_j}s_{i_j}}\text{vol }(\sK_vt_{{\bf a}}^{\bf m}\sK_v)
\varphi_{v, \pi}(t_{\bf a}^{\bf m}).
\end{equation*}
Here $t_{\bf a}^{\bf m}$ corresponds to the vector with
coordinates $a_{i_j}$ at $i_j$. Next

\begin{align*}
\vert b(m)_v \vert & \leq \sum_{1 \leq i_1 < i_2 < \dots i_m \leq
l} \sum_{a_{i_1} >0} \sum_{a_{i_2}>0} \dots \sum_{a_{i_m}>0}
q_v^{-\sum_{j=1}^m a_{i_j}(\sigma_{i_j}-\kappa_{\alpha_{i_j}})} \\
& \leq \sum_{1 \leq i_1 < i_2 < \dots i_m \leq l} \prod_{j=1}^m {
q_v^{ -\sigma_{i_j} + \kappa_{\alpha_{i_j}}}
\over 1-q_v^{ -\sigma_{i_j} + \kappa_{\alpha_{i_j}}}} \\
& \ll \sum_{1 \leq i_1 < i_2 < \dots i_m \leq l} \prod_{j=1}^m
q_v^{ -\sigma_{i_j} + \kappa_{\alpha_{i_j}}} \\
& < \sum_{1 \leq i_1 < i_2 < \dots i_m \leq l} q_v^{m(-1+\epsilon)} \\
& \ll q_v^{m(-1+\epsilon)}.
\end{align*}
Since $m \geq 2$ and $\epsilon < {1 \over 2}$, the series $\sum_v
q_v^{m(1-\epsilon)}$ is convergent, and we are done.

The holomorphy of \eqref{non-triv-eq} follows from the uniform
convergence of
\begin{equation*}
\Sigma_S = \sum_{v \notin S} \vert a_v({\bf s}) \vert
\end{equation*}
on compact subsets of $\cT$, which we now establish. Fix an $\eps'
>0$, and suppose that ${\bf s}\in \cT_{-c+\eps'}$. Then
\begin{equation*}\begin{split}
\Sigma_S & \leq \sum_{v \notin S} \sum_\vartheta
q_v^{-\sigma_\vartheta} \text{vol }(\sK_v
\check\vartheta(\varpi_v) \sK_v) \vert \varphi_v
(\check\vartheta(\varpi_v))) \vert \\
& \ll \sum_{v \notin S} \sum_\vartheta q_v^{-\sigma_\vartheta-c
} \delta (\check\vartheta(\varpi_v)),
\end{split}\end{equation*}
by Lemma of ~\ref{sect:vol} and Theorem~\ref{bound}. Finally,
since by definition
\begin{equation*}
\delta (\check\vartheta(\varpi_v))) = q_v^{\kappa_\vartheta
\ell(\vartheta)},
\end{equation*}
we conclude that
\begin{equation*}
\Sigma_S \ll \sum_{v \notin S} \sum_\vartheta q_v^{-\eta -1} <
\infty,
\end{equation*}
for some $\eta > 0$. This last inequality completes the proof of
the theorem.
\end{proof}

\begin{coro}[of the proof]\label{coro-of-the-proof}
For all $\eps>0$ and all compacts
$K\subset\cT_{-c+\eps}$ there exists a constant
$C(\eps, K)$ such that
\begin{equation*}
\vert \mathcal{J}_S({\bf s}, \pi) \vert \leq C(\eps,K).
\end{equation*}
for all $\pi$ as above, and all ${\bf s} \in K$.
\end{coro}

\begin{coro}\label{coro-inf-int} Let $\K$ be as in Corollary
\ref{coro-spherical-satake}. Let $\phi$ be an automorphic form in
the space of an automorphic representation $\pi$ which is right
invariant under $\K$. Set for $\bfs \in \cT_{\gg 0}$
\begin{equation}\label{inf-int}
\mathcal{J}(\bfs, \phi) : = \int_{\G(\A)} H(\bfs, g)^{-1}\phi(g)dg.
\end{equation}
Then $\mathcal{J}(s, \phi)$ has an analytic continuation to a
function which is holomorphic on $\cT_{-c}$. As usual we let
$\Delta$ be chosen as in the proof of Lemma 4.1 of \cite{Arthur},
and suppose $\phi$ is an eigenfunction for $\Delta$. Define
$\Lambda(\phi)$ by $\Delta\cdot\phi = \Lambda(\phi)\cdot\phi$.
Then for each integer $k>0$, and all $\epsilon > 0$, and every
compact subset $K \subset \cT_{-c+\eps}$, there exists a constant
$C= C(\eps, K, k)$, independent of $\phi$, such that
\begin{equation}
|\mathcal{J}(\bfs, \phi)| \leq C \Lambda(\phi)^{-k} |\phi(e)|
\end{equation}
for all $\bfs \in K$.
\end{coro}
\begin{proof}
Combine Corollary
\ref{coro-of-the-proof}, Corollary \ref{coro-spherical-satake},
Theorem \ref{ramified} and its corollary.
\end{proof}

\begin{rem}
The value of the spherical function at the element
$\check\vartheta(\varpi_v)$, for $\vartheta = r_v
(\iota^*(\alpha))$ is related to the trace of the dominant weight
$\omega_\alpha$ of the $L$-group applied to the Langlands class of
$\pi_v$, and that the Euler product $\mathcal{J}_S({\bf s}, \pi)$ is
regularized by a product of $L_S(s_\alpha, \pi, \omega_\alpha)$.
Additional information about these $L$-functions would lead to
better error terms in asymptotics of rational points.
\end{rem}

\section{Height zeta function}\label{sect:heightzeta}

\subsection{The zeta function} The main tool in the study of distribution properties of
rational points is the height zeta function, defined on $\G (F)
\backslash \G (\A)$ by
$$
\zZ({\bf s}, g) := \sum_{\gamma \in \G(F)} H({\bf s},\gamma g)^{-1}
$$

\begin{prop}
\label{prop:abs-conv} The series defining $\zZ({\bf s},g)$ converges
absolutely to a holomorphic function for ${\bf s} \in \cT_{\gg 0}$.
In its region of convergence
$$
\zZ({\bf s},g) \in C^\infty (\G(F) \backslash \G(\A)).
$$
Furthermore, $\zZ({\bf s},g)$ and all of its group derivatives are
in $\sL^2$.
\end{prop}

\begin{proof}
It suffices to prove the absolute convergence of
$\zZ({\bf s},g)$ for $\Re(s)$ contained in some open cone.
This is  a general fact (see Proposition 4.4 in
\cite{chambert-t02}): since $X$ is projective the cone generated by
ample classes is open in $\Pic(X)_{\R}$. Fix some ample classes
$L_j$ generating $\Pic(X)$. The restriction $\zZ(sL_j,1)$ converges
for $\Re(s)>\sigma_j$, for some $\sigma_j>0$. Now use the
exponential property of heights. The proof of the last statement
is identical to the proof of Proposition 2.3. of \cite{PGL2}.
\end{proof}
\begin{prop}
The function $F(g) = \zZ({\bf s}, g)$ satisfies the conditions of
Lemma \ref{lem:spectral-expansion-general}, and hence has a spectral
expansion.
\end{prop}
\begin{proof} Obvious from the proof of Theorem \ref{thm:one-dim}.
\end{proof}

By Lemma~\ref{lem:spectral-expansion-general}, the zeta function has an
expansion of the form
\begin{equation}
\zZ({\mathbf s}, g) = S_1(\zZ({\mathbf s},\cdot), g) + S_2(\zZ({\mathbf
s},\cdot), g).
\end{equation}
Since the two sides are continuous functions of $g$, we may set
$g=e$ to get
\begin{equation}
\zZ({\mathbf s}) = S_1(\zZ({\mathbf s},\cdot), e) + S_2(\zZ({\mathbf
s},\cdot), e).
\end{equation}
We use this expansion to determine the analytic behavior of the
height zeta function. The idea is to separate out the contribution
of one-dimensional representations.
Since the definition of $S_1$ involves the Laplace operator, the
contribution of one-dimensional representations to this term will
cancel (we will see momentarily that because of the uniform
convergence of the inner sum the $\Delta^n$ can be moved to
the Eisenstein series). It remains to treat the
contribution to $S_2$. Since both sums are of the shape
considered in \eqref{spectral-identity-modified},
without restriction of generality,
we set
\begin{equation}\begin{split}\label{zeta-identity-modified2}
S^\flat({\mathbf s}) = \sum^\flat_{X \in \mX}\sum_{\P} & n(\sA)^{-1} \\
& \int_{\Pi(\M)}\int_{\G(\A)}\left(\sum_{\phi \in \cB_\P(\pi)_X} E(x,
\phi) \overline{E(y, \phi)}\right)H({\mathbf s},y )^{-1} \, dy\,
d\pi,
\end{split}\end{equation}
where the symbol $\flat$ indicates that the summation is over those
classes which do not correspond to one-dimensional representations.
The fact that the innermost sum is uniformly convergent for $y$ in
compact sets is included in the first half of the proof of Lemma 4.4
of \cite{Arthur}. (Note that here too one needs to use Lemma 4.1 of
\cite{Arthur}). Therefore, we may interchange the innermost
summation with the integral over $\G(\A)$ to obtain
\begin{equation}\begin{split}
S^\flat({\mathbf s}) = \sum^\flat_{X \in \mX}\sum_{\P} & n(\sA)^{-1} \\
& \int_{\Pi(\M)}\left(\sum_{\phi \in \cB_\P(\pi)_X} E(e, \phi)
\int_{\G(\A)}\overline{E(y, \phi)}H({\mathbf s},y )^{-1} \,
dy\right) \, d\pi,
\end{split}\end{equation}

\begin{thm}\label{important-theorem}
The function $S^\flat$ has an analytic continuation to a function
which is holomorphic on $\cT_{-c}$.
\end{thm}

\begin{proof} For simplicity we assume that the height function
is invariant under  right and left translation by the compact
subgroup $\K$ as in Corollary \ref{coro-spherical-satake}.
Keeping notation of Corollary \ref{coro-inf-int}, we
write
\begin{equation}
S^\flat({\mathbf s}) = \sum^\flat_{X \in \mX}\sum_{\P}  n(\sA)^{-1}
\int_{\Pi(\M)}\left(\sum_{\phi \in \cB_\P(\pi)_X} E(e, \phi)
{\mathcal J}(\bfs, E(\phi, \cdot))\right) \, d\pi.
\end{equation} We now
use the analytic continuation and bounds established in Corollary
\ref{coro-inf-int} to obtain the analytic continuation of $S^\flat$.
Let $K$ be a compact subset of $\cT_{-c}$. Then there is an
$\epsilon > 0$ such that $K \subset \cT_{-c+\eps}$.
By Corollary \ref{coro-inf-int}, we know that for $\bfs \in K$ and
all $k$ the expression
\begin{equation}
\sum^\flat_{X \in \mX}\sum_{\P}  n(\sA)^{-1}
\int_{\Pi(\M)}\left(\sum_{\phi \in \cB_\P(\pi)_X} |E(e, \phi)|\cdot|
{\mathcal J}(\bfs, E(\phi, \cdot))|\right) \, d\pi
\end{equation}
is bounded by
\begin{equation}
C(\eps, K, k) \sum^\flat_{X \in \mX}\sum_{\P}  n(\sA)^{-1}
\int_{\Pi(\M)}\left(\sum_{\phi \in \cB_\P(\pi)_X} \Lambda(\phi)^{-k}
|E(e, \phi)|^2 \right) \, d\pi.
\end{equation}
The convergence of the last expression is a consequence of
Proposition \ref{important-proposition}. This establishes Theorem
\ref{important-theorem}. \end{proof}
\begin{rem}
In the anisotropic situation, the desired analytic properties follow
from the analytic continuation of the spectral zeta function of the
Laplace operator on the corresponding compact quotient.
\end{rem}
\subsection{}
Now we let $\underline{a} = (a_\alpha) \in \mathcal{T}_\N$, and
for $s \in \C$, we set
\begin{equation}
\zZ_{\underline{a}}(s) = \zZ(s(a_1, a_2, \dots, a_r))
\end{equation}
as a function of one complex variable. We need to determine the
right most pole of $\zZ_{\underline{a}}(s)$. Set
\begin{equation}
\sigma(\underline{a}) = \max_{i} \frac{1 +
\kappa_{\alpha}}{a_\alpha},
\end{equation}
and let $S(\underline{a})$ be the set of $\alpha$, modulo Galois
action, for which the maximum is achieved, and $m(\underline{a})=
\# S(\underline{a})$. The theorem implies that
$\zZ_{\underline{a}}(s)$ has no pole for $\Re(s) >
\sigma(\underline{a})$. The order of pole of
\begin{equation}
\int_{\G(\A)}H(s \underline{a}, g)^{-1}\, dg
\end{equation}
at $s = \sigma(\underline{a})$ is equal to $m(\underline{a})$. Therefore,
we need those automorphic characters $\chi$ such that
$\xi_{\alpha}(\chi) =1 $ for all $\alpha \in S(\underline{a})$.
Clearly, we are interested only in those automorphic
$\chi$ which satisfy
\begin{equation}
\int_{\G(\A)} H(s, g)^{-1} \chi(g) \, dg \ne 0
\end{equation}
for some $s$ in the domain of absolute convergence. This implies
that $\chi$ is right, and in this case also
left, invariant under the compact open subgroup $\sK$ of $\G(\A_f)$.
Let $X(\underline{a})$ be the collection of all such characters.
By Lemma 3.1 of \cite{GMO}, the set $X(\underline{a})$ is
finite. The proof of the above theorem shows that
\begin{thm}\label{one-dim-limit}
The complex function $\zZ_{\underline{a}}(s)$ has a meromorphic
continuation to $\Re(s)>\sigma(\underline{a})- \epsilon$,
$\epsilon
> 0$, with an isolated pole at $\sigma(\underline{a})$ of order $m(\underline{a})$.
Furthermore,
\begin{equation*}
\lim_{s \to \sigma(\underline{a})}
(s-\sigma(\underline{a}))^{m(\underline{a})} \zZ_{\underline{a}}(s)
=  \lim_{s \to \sigma(\underline{a})}
(s-\sigma(\underline{a}))^{m(\underline{a})}\sum_{\chi \in
X(\underline{a})}\int_{\G(\A)} H(s, g)^{-1}\chi(g) \, dg.
\end{equation*}
The limit is a positive real number.
\end{thm}
\begin{proof} If we set
\begin{equation}
\G_{\underline{a}} = \bigcap_{\chi \in X(\underline{a})} \ker(\chi),
\end{equation}
then the above limit is equal to
\begin{equation}
\lim_{s \to \sigma(\underline{a})}
(s-\sigma(\underline{a}))^{m(\underline{a})}\int_{\G_{\underline{a}}}
H(s, g)^{-1} \, dg.
\end{equation}
Again by Lemma 3.1. of \cite{GMO}, we know that the group
$\G_{\underline{a}}$ has finite index in $\G(\A)$. This easily
implies that the above limit is a positive number.
\end{proof}
\begin{rem}
This result is a fundamental ingredient in the treatment of
\cite{GMO}.
\end{rem}

A special case of particular interest is when
$a_\alpha = \kappa_{\alpha} + 1$, for all $\alpha$. Because of its
relevance to the anticanonical class of the wonderful
compactification (c.f. Proposition \ref{prop:geometry}), we denote
this $\underline{a}$ by $\underline{\kappa}$. In this case,
$\sigma(\underline{\kappa}) = 1$, and $m(\underline{\kappa})$ is
equal to the number of distinct Galois orbits in $\Delta(\G, \T)$.

\begin{prop}\label{pole-hecke}
Suppose $\chi$ is an automorphic character such
that $\xi_\alpha(\chi) =1$ for all $\alpha$. Then $\chi \equiv 1$.
\end{prop}

\begin{proof} This easily follows from Proposition~\ref{prop:41.1},
Cartan decomposition, and the weak approximation.
\end{proof}

The proposition shows that $X(\underline{\kappa}) = \{ 1 \}$.

\begin{thm}\label{thmmain2}
The complex function $\zZ_{\underline{\kappa}}(s)$
has a meromorphic continuation to $\Re(s) >
1 - \epsilon$ for some $\epsilon > 0$ with an isolated pole of order
equal to the number of Galois orbits in $\Delta(\G, \T)$ at $s=1$.
Furthermore,
\begin{equation*}
\lim_{s \to \sigma(\underline{\kappa})}
(s-\sigma(\underline{\kappa}))^{m(\underline{\kappa})}
\zZ_{\underline{\kappa}}(s) =  \lim_{s \to \sigma(\underline{\kappa})}
(s-\sigma(\underline{\kappa}))^{m(\underline{\kappa})}
\int_{\G(\A)} H(s\underline{\kappa},g)^{-1} \, dg.
\end{equation*}
\end{thm}

\subsection{Examples}

\begin{exam}[$\PGL_n$]\label{ex:pgln}
Let $F_v$ be a local field, $\G=\PGL_n$ and $\chi_v$ a
one-dimensional representation of $\G(F_v)$. Then $\chi_v$ has the
form
$$
\chi_v(g)=\xi_v(\det(g)),
$$
where $\xi_v$ is a character of $F^*_v$ whose order divides $n$.
In the global situation,
$$
\chi(g)=\xi(\det(g)),\,\,\, g\in \G(\A),
$$
with $\xi$ of order dividing $n$, and $\xi|_{F^*} =1$.
Let $\alpha_1, \dots, \alpha_{n-1}$ be the simple roots, with
the convention that
\begin{equation}
\alpha_i (\text{diag }(a_1, \dots, a_n) ) = a_i / a_{i+1},
\end{equation}
and
\begin{equation}
\check\alpha_i(t) = \begin{pmatrix} t I_i \\ & I_{n-i}.
\end{pmatrix}
\end{equation}

It is not hard to see that
\begin{equation}
2 \rho = \sum_{i=1}^{n-1} i(n-i) \alpha_i.
\end{equation}
Let $\underline{a} = (a_1, \dots, a_{n-1}) \in \N^{n-1}$. We
define $\sigma( \underline{a})$, $S(\underline{a})$,
$m(\underline{a})$, and $X(\underline{a})$ as above. Suppose $\chi
\in X(\underline{a})$. Then if $\chi = \xi \circ \det$, we must
have
\begin{equation}
\xi^i = 1
\end{equation}
for all $i \in S(\underline{a})$. In particular, if we set
\begin{equation}
d(\underline{a}) = \gcd(n, \gcd_{i \in S(\underline{a})} i),
\end{equation}
then we obtain that $\chi \in X(\underline{a})$, if and only if
$\chi^{d(\underline{a})} = 1$.
\end{exam}

\begin{exam}[$\PGL_4$]
Here
$$
2\rho (\diag(a,b,c,d)) =c^3bc^{-1}d^{-3} = \al^3\beta^4\gamma^3,
$$
if $\Delta=\{ \al=a/b,\beta=b/c,\gamma=c/d\}$. Thus
$\underline{a}_{-K} = (4,5,4)$. For $\underline{a}=(1,1,1)$
we have
$$
\sigma(\underline{a})=\max(4,5,4)=5 \,\,\, \text{ and } \,\,\,
S(\underline{a}) =\{ \beta \}.
$$
We have
$$
\check{\beta}(\varpi) =\diag(\varpi,\varpi,1,1)\,\,\,\text{ and }
\,\,\, \chi_v(\diag(\varpi,\varpi,1,1)) =\xi_v(\varpi^2).
$$
Let $\xi$ have order two. By our analysis, $\chi = \xi \circ
\det$ contributes to the asymptotic constant.
\end{exam}

\begin{exam}[${\rm PGU}_3$]
Let $E/F$ be a quadratic extension and ${\rm GU}_3\subset
\GL_3(E)$ the set of all $g$ such that
$$
{}^{T}g^{\sigma} S g=\la S,\,\,\, \text{ where } \,\,\,
S=\left(\begin{array}{ccc} & & 1  \\ & 1 & \\ 1 & &
\end{array}\right).
$$
Let $\G={\rm PGU}_3$ and write $\la=\nu(g)$ for the similitude
norm of $g\in \G$. Then
$$
\G^{sc}={\rm SU}_3, \sZ=\sZ^{sc}=\{ \zeta \cdot I_3,\,\,  \zeta^3=1\}.
$$
The Galois action on $\G^{sc}$ is given by
$$
\tilde{\rho}(g)=S {}^{T}g\rho S^{-1}
$$
if $\rho\in \Ga_{\bar{F}/F}$ restricts to $\sigma \in \Ga_{E/F}$
and $\tilde{\rho}(g)=g^{\rho} $ if  $\rho\in \Ga_{\bar{F}/E}$.

\

We claim that $j\,:\, \T^{sc}(F_v)\ra \T(F_v)$ is not surjective.

\

We have
$$
\T^{sc}(F_v)=\{ \diag(a, \eps, a) \,\,\, a\in E_w^*\}
$$
and $aa^{-\sigma}\eps =1$, with $\Nm(\eps)=1$, whereas
$$
\T(F_v)=\{ \diag(a, b, \la a^{-\sigma}), \,\,  \text{ modulo
scalars } \}
$$
with $\la=b^{\sigma} b$. We can take $t\in \T(F_v)$ to have the
form
$$
\diag(a, 1, a^{-\sigma}).
$$
Finally,
$$
j(\diag(a, 1, a^{-\sigma}))= \diag(a\eps^{-1}, 1,
a^{-\sigma}\eps^{\sigma}) = \diag(u, 1, u^{-\sigma}),
$$
modulo scalars. One-dimensional automorphic representations of
$\G(\A)$ are given by
$$
\chi(g)=\xi(\nu(g)).
$$
We need to consider the range of $\nu$.

If $v$ remains prime in $E$, then
$$
{}^{T}g_v^{\sigma} Sg_v =\la_v\cdot S
$$
which implies $\la_v^3=\Nm_{w/v}(\det(g_v))$. Thus $\la_v\in
\Nm_{w/v}(E^*_w)$. If $v$ splits in $E$, then there is no
condition on $\la_v$. Thus
$$
\nu\,:\, \G(\A)\ra \Nm_{E/F} (\A_E^*).
$$
Also $\nu(\al_v\cdot I_3) = \Nm_{W/V} (\al_v) $ if $v$ remains
prime in $E$.

We may take for the adjoint group ${\rm U}_3 /\sZ$, where
$\sZ=\diag(\eps,\eps,\eps)$ and $\Nm_{E/F}(\eps)=1$. We have
$\chi(g)=\xi(\det(g))$, $g\in \G(\A)$. Here $\xi\in
\N^1(\A)/\N^1(F)$, $\xi^3=1$. By Hasse's theorem,
$$
\begin{array}{ccc}
\mathbb G_m(\A_E) & \ra &  \N^1(\A) \\
    \mathbb G_m(E) & \ra & \N^1(E) \\
     x & \mapsto & x\cdot x^{-\sigma}
\end{array}
$$
are surjective. We may view $\xi $ as a character $\eta$ of
$\mathbb G_m(\A_E)/\mathbb G_m(E)$ with $\eta^3=1$. The map
$$
j\,:\, \T^{sc}(F_v)\ra \T(F_v)
$$
is not surjective, since $\T(F_v)$ is the set of all
$\diag(a,1,a^{-\sigma})$ and $\T^{sc}(F_v)$ the set of
$\diag(a,\eps,a^{-\sigma})$, with $aa^{-\sigma}\eps =1$.

We wish to define
$$
L(s,\chi)=L_E(s,\eta).
$$
This leads to the necessity of introducing a norm map. In this
example we have
$$
\begin{array}{rcl}
\Nm_{w/v} \,:\, \T(E_w) & \ra &  \T(F_v) \\
                   t   & \mapsto  & t\cdot \tilde{\sigma}(t)
\end{array}
$$
with $\tilde{\sigma}(t) = St^{\sigma} S^{-1}$. Thus
$$
\Nm_{w/v} (\diag(a,b,c) = \diag(ac^{-\sigma}, bb^{-\sigma},
ca^{-\sigma}).
$$
Note that
$$
1\ra {\rm SU}_v/\mu_3 \ra {\rm U}_3/\sZ\stackrel{\det}{\lra}
\N^1_v/(\N^1_v)^3 \ra 1
$$
is exact.
\end{exam}

\section{Manin's conjecture}

In this section we apply Theorem~\ref{thmmain2} to obtain a
result regarding the number of rational points of bounded height
on wonderful compactifications of the group $\G$.

We first recall a standard Tauberian theorem.

\begin{thm}
\label{tauberian-thm} Let $(c_n)_{n\in \N}$ be a sequence of
positive real numbers such that for any $B>0$ only finitely many
$c_n<B$. Let $f$ be a function with values in $\R_{>0}$ such that
\begin{equation*}
Z(s,f) = \sum_{n} \frac{f(n)}{c_n^s}
\end{equation*}
is absolutely convergent for $\Re(s) > a> 0$ and has a
representation
\begin{equation*}
Z(s,f) = {1 \over (s-a)^b}g(s) + h(s)
\end{equation*}
with $g(s)$ and $h(s)$ holomorphic for $\Re(s) \geq a$, $g(a)
\ne0$ and $b \in \N$. Then
\begin{equation*}
\sum_{n\,:\, c_n\leq B} f(n) = {g(a) \over a(b-1)!}B^a (\log
B)^{b-1}(1 + o(1)),
\end{equation*}
as $B \to \infty$.
\end{thm}

\begin{thm}
\label{thm:main2} Let $X$ be the compactification of a semi-simple
group $\G$ of adjoint type over $F$ as in Section~\ref{sect:geometry}
and $\cL=(L,\|\cdot \|_v)$ an adelically metrized line bundle
such that its class $[L]\in\Pic(X)$ is contained in the interior of
the cone of effective divisors $\La_{\rm eff}(X)$. Then
$$
N(\G,\cL,B):=  \# \{ x\in \G(F)\,|\, H_{{\cL}}(x)\le B\} =
\Theta(\cL) B^{a(L)} \log(B)^{b(L)-1}(1+o(1))
$$
as $B\ra \infty$. Here
$$
a(L)=\inf\{ a\,\, | \,\, a[L]+[K_X]\in \La_{\rm eff}(X)\}
$$
(where $K_X$ is the canonical line bundle of $X$) and $b(L)$ is
the (maximal) codimension of the face of $\La_{\rm eff}(X)$
containing $a(L)[L]+[K_X]$. Moreover, $\Theta(-\cK_X)$ is the
constant defined in \cite{peyre95}.
\end{thm}

\begin{proof}
We combine the following facts:
\begin{itemize}
\item $\Br(X)/\Br(F)=1$;
\item $\overline{G(F)}=G(\A)$;
\item $\int_{G(\A)} \omega_S = \int_{X(\A)}\omega_S$ (where $\omega_S$
is the Tamagawa measure defined in \cite{peyre95};
\item $\omega_v
= \tau_v(G)^{-1} H_{-K_X}(\cdot)dg_v$, for all $v\in \Val(F)$.
\end{itemize}
The first two facts are in \cite{sansuc}, here it is important that
$\G$ is of adjoint type.
For the identifications of the measures see \cite{peyre95},
\cite{BT}, Theorem 7.2 or \cite{batyrev-t98}.
\end{proof}

\end{document}